\documentclass[a4paper]{amsart}
\usepackage{mathrsfs,amssymb}
\usepackage{multicol}\multicolsep=0pt

\def\crulefill{\leavevmode\leaders\hrule height 1pt\hfill\kern 0pt}
\long\def\QUERY#1{%
\leavevmode\newline%
\noindent$\star\star\star$\thinspace\textsf{Comment/Query}\crulefill\newline%
   \space #1\newline\hbox to 120mm{\crulefill}$\star\star\star$\newline
}

\numberwithin{equation}{section} \theoremstyle{definition}
\newtheorem{Defn}[equation]{Definition}

\theoremstyle{plain}
\newtheorem{Prop}[equation]{Proposition}
\newtheorem{Theorem}[equation]{Theorem}

\newtheorem{Lemma}[equation]{Lemma}
\newtheorem{Cor}[equation]{Corollary}

\def\Case#1{\medskip\noindent\textbf{Case #1}:\leavevmode\newline}
\def\subcase#1{\bigskip\noindent\textbf{Subcase #1}:\leavevmode\newline}

\makeatletter
\def\enumerate{\begingroup\ifnum\@enumdepth>3\@toodeep\else
      \advance\@enumdepth\@ne
      \edef\@enumctr{enum\romannumeral\the\@enumdepth}%
      \topsep\z@\parskip\z@
      \list{\csname label\@enumctr\endcsname}
        {\@nmbrlisttrue\let\@listctr\@enumctr
         \parsep\z@\itemsep\z@\topsep\z@
         \setcounter{\@enumctr}{0}
         \def\makelabel##1{\hss\llap{\rm ##1}}
       }\fi}

\makeatother

\let\bar=\overline
\let\epsilon=\varepsilon
\def\({\big(}
\def\){\big)}

\def\R{\mathbb R}

\def\M{\mathfrak M}

\def\0{\underline{0}}

\def\H{\mathscr H}

\def\Sym{\mathfrak S}

\DeclareMathOperator{\Char}{char} 
\DeclareMathOperator{\rank}{rank}

\def\simk{\overset{k}\sim}

\def\simnn{\overset{n-1}\sim}

\def\floor#1{\lfloor\tfrac#1\rfloor}
\def\UPD{\mathscr{T}^{ud}}
\def\Std{\mathscr{T}^{std}}

\def\s{\mathfrak s}
\def\ts{\tilde\s}
\def\t{\mathfrak t}
\def\u{\mathfrak u}
\def\v{\mathfrak v}
\def\q{\mathbf q}
\def\r{\mathbf r}

\def\ss{\mathbf s}
\def\ts{\mathbf t}
\def\us{\mathbf u}

\newcommand{\cba}[1]{\mathscr {B}_{#1}}

{\catcode`\|=\active
  \gdef\set#1{\mathinner{\lbrace\,{\mathcode`\|"8000%
                                   \let|\midvert #1}\,\rbrace}}
  \gdef\seT#1{\mathinner{\Big\lbrace\,{\mathcode`\|"8000%
                                   \let|\midverT #1}\,\Big\rbrace}}
}
\def\midvert{\egroup\mid\bgroup}
\def\midverT{\egroup\,\Big|\,\bgroup}

\def\Set[#1]#2|#3|{\Big\{\ #2\ \Big| \
           \vcenter{\hsize #1mm\centering #3}\Big\}}

\title{Gram determinants and semisimplicity  criteria for Birman-Wenzl algebras}
\author{Hebing Rui and Mei Si}
\date{Revised version, Feb. 18, 2008}
\address{H.R. Department of Mathematics,  East China Normal
University, Shanghai, 200062, China}
\email{hbrui@math.ecnu.edu.cn}
\address{M.S. Department of Mathematics,  East China Normal
University, Shanghai, 200062, China}
\email{52050601011@student.ecnu.edu.cn}

\thanks{The first author  is  supported in part by NSFC
and NCET-05-0423}

\begin{document}
\baselineskip15pt
\begin{abstract} In this paper, we compute all Gram determinants
associated to all cell modules of  Birman-Wenzl algebras. As a
by-product, we give a  necessary and sufficient condition for
Birman-Wenzl algebras being semisimple over an arbitrary field.
\end{abstract}

\sloppy \maketitle

\section{Introduction}

In \cite{BM}, Birman and Wenzl introduced a class of associative
algebras $\cba{n} $, called Birman-Wenzl algebras, in order to study
link invariants. They are  quotient algebras of the group algebras
of  braid groups. On the other hand, there is a Schur-Weyl duality
between $\cba{n}$ with some special parameters over  $\mathbb C$ and
quantum groups of types $B, C, D$~\cite{W2}. Thus, $\cba{n}$ plays
an important role in different disciplines.

In this paper, we  work on $\cba{n}$ over the ground ring
$R:=\mathbb Z[r^\pm, q^\pm, \omega^{-1}]$ where $\omega=q-q^{-1}$
and $q, r$  are indeterminates.

\begin{Defn}\cite{BM}\label{bmw-def} The Birman-Wenzl algebra
$\cba{n}$ is a unital associative $R$-algebra with generators
$T_i,1\le i\le n-1$ and relations\begin{enumerate} \item
$(T_i-q)(T_i+q^{-1})(T_i-r^{-1})=0$, for $1\le i\le n-1$,
\item $T_{i}T_{i+1}T_i=T_{i+1}T_iT_{i+1}$,  for $1\le i\le n-2$,
\item $ T_iT_j=T_jT_i$,  for $|i-j|>1$,
\item $E_iT_j^\pm E_i=r^\pm E_i$,  for $1\le i\le n-1$ and $j=i\pm 1$,
\item $ E_i T_i=T_iE_i=r^{-1} E_i$,  for $1\le i\le n-1$,
\end{enumerate}
where $ E_i=1-\omega^{-1}(T_i-T_{i}^{-1})$ for $1\le i\le n-1$.
\end{Defn}

In \cite{MW}, Morton and Wassermann proved that $\cba{n}$ is
isomorphic to  Kauffman's tangle algebra~\cite{K} whose $R$-basis is
indexed by Brauer diagrams. This enables them to show that $\cba{n}$
is a free $R$-module with rank $(2n-1)!!$. Let $F$ be a field which
contains non-zero elements $\mathbf q, \mathbf r$ and $\mathbf
q-\mathbf q^{-1}$. Then the  Birman-Wenzl algebra $\mathscr B_{n,
F}$ over $F$ is isomorphic to $\cba{n}\otimes_{R} F$. In this case,
$F$ is considered as an $R$-module such that $r, q, \omega$ act on
$F$ as $ \mathbf r$, $\mathbf q$, and $\mathbf q-\mathbf q^{-1}$,
respectively. We will use $\cba{n}$ instead of $\mathscr B_{n, F}$
if there is no confusion.

Let $\langle E_1\rangle$ be the two-sided ideal of $\cba{n}$
generated by $E_1$. It is well-known that  $\cba{n}/\langle
E_1\rangle$ is isomorphic to the Hecke algebra $\mathscr H_n$
associated to the symmetric group $\mathfrak S_n$. If we denote by
$g_i, 1\le i\le n-1$ the distinguished generators of $\mathscr H_n$,
then the defining relations for  $\H_n $ are as follows:
$$
\begin{aligned}(g_i-q)(g_i+q^{-1})&=0 \text{ for $1\le i\le
n-1$,}\\
 g_{i}g_{i+1}g_i&=g_{i+1}g_ig_{i+1},  \text{ for $1\le i\le
 n-2$,}\\
 g_ig_j&=g_jg_i,  \text{ for $|i-j|>1$.}\end{aligned}
$$
The corresponding isomorphism from $\cba{n}/\langle E_1\rangle$ to
$\mathscr H_n$ sends $T_i\pmod {\langle E_1\rangle}$ to $g_i$ for
all $1\le i\le n-1$.

Using  the cellular structure of $\mathscr H_n$ together with
Morton-Wassermann's result in \cite{MW}, Xi~\cite{Xi} proved that
$\cba{n}$ is  cellular  over $R$ in the sense of \cite{GL}.

In \cite{GL}, Graham and Lehrer constructed a class of generically
irreducible modules for each cellular algebra, which are  called
cell modules. A question arises. When is a generically  irreducible
cell module not irreducible? Graham and Lehrer proved that  any cell
module of a cellular algebra is equal to its simple head if  and
only if the cellular algebra is (split) semisimple. This gives a
method to determine the semisimplicity of a cellular algebra.

There is no result on the first problem for $\cba{n}$. In \cite{W2},
Wenzl used the  ``Jones basic construction'' and the Markov trace on
$\cba{n}$ to give   some partial results for  $\cba{n}$ being
semisimple. More explicitly, Wenzl~\cite[5.6]{W2} proved that
 $\cba{n}$  is semisimple over $\mathbb C$ except possibly if $q$ is a root
 of unity or $r=q^k$ for some $k\in \mathbb Z$. However, there is no explicit
description for such  $k$'s.

Enyang constructed  the Murphy basis for each cell module of
$\cba{n}$ in \cite{Enyang}  on which the Jucys-Murphy elements of
$\cba{n}$ act upper triangularly. This enables us to use standard
arguments (see e.g. \cite{JM} or more generally, \cite{M:semi})  to
construct an  orthogonal basis of $\cba{n}$. Via this orthogonal
basis together with   classical branching rule for $\cba{n}$ in
\cite{W2}, we obtain a recursive formula for the Gram determinant
associated to each cell module of $\cba{n}$. This is the first main
result of this paper.

Let $\Lambda^+(n)$ be the set of all partitions of $n$.  When
$r\not\in \{q^{-1}, -q\}$,  we will prove that $\cba{n}$ is
semisimple if and only if $$\prod_{k=2}^n \det G_{1, (k-2)} \det
G_{1, (1^{k-2})} \prod_{\lambda\in \Lambda^+(n)}\det G_{0,
\lambda}\neq 0 .$$ Using our recursive formulae on Gram
determinants, we compute $\det G_{1, \lambda}$ explicitly for
$\lambda\in \cup_{k=2}^n \{(k-2), (1^{k-2})\}$. Note that
$\prod_{\lambda\in \Lambda^+(n)} \det G_{0, \lambda}\neq 0$  if and
only if $\H_n$ is semisimple. So, we can give a criterion for
$\cba{n}$ being semsimple when $r\not\in \{q^{-1}, -q\}$. When $r\in
\{q^{-1}, -q\}$, we can determine whether $\cba{n}$ is semisimple by
elementary computation. It gives a complete solution of the problem
on the semisimplicity of $\cba{n}$ over an arbitrary field. This is
the second main result of the paper.

Note that the group algebra of  $\mathfrak S_n$ is both a subalgebra
and a quotient algebra of the Brauer algebra $B_n$ \cite{BrauerAlg}.
Thus, Doran-Wales-Hanlon\cite{DoranHanlonWales} can restrict a
module for
 $B_n$ to the group algebra of $\mathfrak S_n$. However,
 $\mathscr H_n$ is not a subalgebra of $\cba{n}$. We can not
restrict a $\cba{n}$-module to $\mathscr H_n$. In other words, we
can not use the method in \cite{Rui:ssbrauer, RS:ssbrauer} to give a
criterion for $\cba{n}$ being semisimple. Finally,  we remark that
the method we use in the current paper can be used to deal with
cyclotomic Nazarov-Wenzl algebras~\cite{AMR}. Details will appear
elsewhere.

We organize this paper as follows. In section~2, we recall the
Jucys--Murphy basis for each cell module of $\cba{n}$  in
\cite{Enyang}.  An orthogonal basis of each cell module of $\cba{n}$
will be constructed in section~3. In section~4, we prove the
recursive formulae on Gram determinants. Finally, we give a
criterion for $\cba{n}$ being semisimple in  section~5.

\section{Jucys-Murphy basis for $\cba{n}$}
In this section, unless otherwise stated,  we assume  $R=\mathbb
Z[r^{\pm}, q^{\pm}, \omega^{-1}]$ where $\omega=q-q^{-1}$ and $q, r$
are indeterminates. The main purpose of this section is to construct
the Jucys-Murphy basis of $\cba{n}$ by using Enyang's basis of each
cell module of $\cba{n}$. We state some identities  needed later on.
We start by recalling the definition of Jucys-Murphy elements $L_i,
1\le i\le n$ for $\cba{n}$ in \cite{Enyang}.

Define $L_1=r$ and $L_i=T_{i-1}L_{i-1}T_{i-1}$ for $2\le i\le
n$\footnote{In \cite{Enyang}, Enyang defined $L_1=1$ and
$L_i=T_{i-1}L_{i-1}T_{i-1}$. }. The following identities can be
found in \cite{BM} and \cite{Enyang}.

\begin{Lemma}[\cite{BM,
Enyang}] \label{wenzl-rel} Suppose
$\delta=\frac{(q+r)(qr-1)}{r(q+1)(q-1)}$. We have:
\begin{enumerate}\item
$E_i^2=\delta E_i$, $1\le i\le n-1$,
\item $E_i T_j=T_j E_i$,  $|i-j|>1$,
\item $T_i^2=1+\omega (T_i-r^{-1}E_{i})$, $1\le i\le n-1$,
\item
$E_iE_jE_i=E_i$ for  $1\le i\le n-1$ and $j=i\pm 1$,
\item
$E_i E_{j}=T_{j}T_i E_{j}=E_iT_{j}T_{i}$ for  $1\le i\le n-1$ and
$j=i\pm 1$,
\item $T_i L_k=L_kT_i$ if $k\not\in \{i, i+1\}$,
\item $E_i L_k=L_k E_i$ if $k\not\in \{i, i+1\}$,
\item $L_iL_k=L_kL_i$ for all $1\le i, k\le n$,
\item $T_i L_iL_{i+1}=L_{i} L_{i+1} T_i$ for all $1\le i\le n-1$,
\item $L_2L_3\cdots L_n$ is a central element in $\cba{n}$.
\end{enumerate}
\end{Lemma}

The following result is well-known. One can prove it by  checking
the defining relations for $\cba{n}$ in Definition~\ref{bmw-def}.

\begin{Lemma}\label{anti}  \begin{enumerate}\item There is a quasi
$R$-linear automorphism  $\sigma: \cba{n}\rightarrow \cba{n}$ such
that $\sigma(T_i)=T_{i}^{-1}$, $\sigma(q)=q^{-1}$,
$\sigma(r)=r^{-1}$. Therefore, $\sigma(\delta)=\delta$,
$\sigma(E_i)=E_i$ and $\sigma(L_j)=L_j^{-1}$ for $1\le i\le n-1$
and $1\le j\le n$.
\item There is an  $R$-linear anti-involution $\ast: \cba{n}\rightarrow
\cba{n}$ such that $T_i^\ast =T_{i}$. Thus, $E_i^\ast=E_i$ and
$L_j^\ast =L_j$ for $1\le i\le n-1$ and $1\le j\le
n$.\end{enumerate}
\end{Lemma}

\begin{Lemma} \label{rs-rel} Suppose that $k$ is a positive integer. The following
equalities hold.
\begin{enumerate}\item $L_iL_{i+1}E_i=E_i=E_{i}L_iL_{i+1}$,
$1\le i\le n-1$.
\item  $T_{i-1} L_i^k =L_{i-1}^k T_{i-1}+\omega\sum_{j=1}^k
L_{i-1}^{j-1}(1-E_{i-1}) L_i^{k-j+1}$, $2\le i\le n$.
\item  $T_{i} L_i^k =L_{i+1}^k T_{i}-\omega\sum_{j=1}^k
L_{i+1}^{j}(1-E_{i}) L_i^{k-j}$, $1\le i\le n-1$.
\item  $T_{i-1} L_i^{-k} =L_{i-1}^{-k} T_{i-1}-\omega\sum_{j=1}^k
L_{i-1}^{-j}(1-E_{i-1}) L_i^{j-k}$,   $2\le i\le n$.
\item  $T_{i} L_i^{-k} =L_{i+1}^{-k} T_{i}+\omega\sum_{j=1}^k
L_{i+1}^{1-j}(1-E_{i}) L_i^{j-k-1}$,   $1\le i\le n-1$.
\end{enumerate}
\end{Lemma}
\begin{proof} (a) can be proved by induction on $i$.  Note that $
 T_{i-1} L_i = T_{i-1}^2 L_{i-1} T_{i-1}$. Now, (b) follows from
 Definition~\ref{bmw-def}(e) and Lemma~\ref{wenzl-rel}(c) for
 $k=1$. In general, by induction assumption,
\begin{equation}\label{tl}
T_{i-1} L_i^k = T_{i-1}L_i^{k-1} L_i=( L_{i-1}^{k-1} T_{i-1} +
\omega \sum_{j=1}^{k-1}L_{i-1}^{j-1}(1-E_{i-1}) L_i^{k-j})L_i.
\end{equation}
 So, (b) follows if we use (b) for
$k=1$ to rewrite $T_{i-1} L_i$ in (\ref{tl}). One can verify (c)
similarly. Applying  $\sigma$ to (b) (resp.  (c)) and using
Lemma~\ref{wenzl-rel}(c) yields (d) (resp. (e)).\end{proof}

\begin{Lemma}\label{e-powerk} For any $1\le i\le n-1$, and $k\in \mathbb N$,
$$
E_iL_i^kE_i=r^2 E_i L_i^{-k} E_i+r\omega\sum_{j=1}^{k-1}
(E_iL_i^{-j}E_iL_i^{k-j} E_i-E_iL_i^{k-2j}E_i)
$$
\end{Lemma}
\begin{proof} By induction on $i$, we have  \begin{equation}\label{e-power1}
E_iL_iE_i=r(\delta+\omega\sum_{j=1}^{i-1}
(L_j-L_j^{-1}))E_i.\end{equation} Applying  $\sigma$ to $E_iL_iE_i$
and using (\ref{e-power1}) yields $r^2E_iL_i^{-1} E_i=E_iL_i E_i$.
This proves  the result for $k=1$. In general,  we have $E_i L_i^k
E_i =rE_i T_iL_i^k E_i$. Using Lemma~\ref{rs-rel}(a)(c) and
Definition~\ref{bmw-def} to simplify $E_i T_iL_i^k E_i$ yields the
formula as required.
\end{proof}

For any $R$-algebra $A$, let $Z(A)$ be the center of $A$.
\begin{Prop}\label{omega}
Given a positive integer   $ i\le n-1$ and  an integer $k$. We have
$E_iL_i^k E_i=\omega_i^{(k)} E_i$, where
 $\omega_i^{(k)}\in R[L_1^\pm, L_2^\pm , \cdots, L_{i-1}^\pm]\cap
Z(\cba{i-1})$.
\end{Prop}
\begin{proof} By Lemma~\ref{e-powerk}, we can assume that $k\ge 0$
without loss of generality.  The case $k=0$ is trivial since
$E_i^2=\delta E_i$. We prove the result by induction on $i$ and $k$
for $k>0$. Since we are assuming that $L_1=r$,
$\omega_1^{(k)}=r^k\delta$. When $k=1$,  the result follows from
(\ref{e-power1}). Now, we assume that $i>1$ and $k>1$.

Write $L_i^k =T_{i-1} L_{i-1} T_{i-1} L_{i}^{k-1}$.  By
Lemma~\ref{rs-rel}(b),
\begin{equation}\label{center111}
E_iL_i^kE_i=E_iT_{i-1}L_{i-1}^k T_{i-1}E_i+\omega\sum_{j=1}^{k-1}
E_{i}T_{i-1}L_{i-1}^j (1-E_{i-1}) L_i^{k-j}E_i
\end{equation}

First, we consider $\sum_{j=1}^{k-1} E_i T_{i-1}L_{i-1}^j
(1-E_{i-1}) L_i^{k-j} E_i$. Applying $\ast$ to Lemma~\ref{rs-rel}(c)
yields $ L_{i-1}^k T_{i-1}=T_{i-1}L_i^k -\omega\sum_{j=1}^k
L_{i-1}^{k-j}(1-E_{i-1}) L_i^j$.
 Multiplying $E_iE_{i-1}$ (resp. $E_i$) on the left (resp. right)
 of $L_{i-1}^k T_{i-1}$ and using Definition~\ref{bmw-def}(e),
 Lemma~\ref{rs-rel}(a) and Lemma~\ref{wenzl-rel}(d) together with induction assumption
 on $E_{i-1}L_{i-1}^{j}E_{i-1}$ for $j<k$,  we have
\begin{equation}\label{center} E_iE_{i-1}L_{i-1}^k T_{i-1}E_i
=r^{-1} L_{i-1}^{-k}E_i-\omega
 \sum_{j=1}^{k}
 (L_{i-1}^{k-2j}-\omega_{i-1}^{(k-j)}L_{i-1}^{-j}) E_{i}\end{equation}

Similarly, we have
\begin{equation}\label{center1} E_iT_{i-1}L_{i}^k E_i
=r L_{i-1}^{k}E_i+\omega
 \sum_{j=1}^{k}
 (L_{i-1}^{j-1}\omega_i^{(k-j+1)}-L_{i-1}^{2j-2-k}) E_{i}\end{equation}

Applying $\ast$ on both
sides of (\ref{center}) and using (\ref{center1}), we have
\begin{equation}\label{f1} \sum_{j=1}^{k-1} E_{i}T_{i-1}L_{i-1}^j
(1-E_{i-1}) L_i^{k-j}E_i=E_i f_1
\end{equation} for some  $f_1\in R[L_1^\pm, L_2^\pm , \cdots,
L_{i-1}^\pm]$. Now, we   discuss $E_iT_{i-1}L_{i-1}^k T_{i-1}E_i$.
We have
$$\begin{aligned} E_iT_{i-1}L_{i-1}^k T_{i-1}E_i& =E_iE_{i-1}T_{i}^{-1} L_{i-1}^k
T_{i-1}E_{i}\\
&=E_iE_{i-1}  (T_i-\omega (1-E_i))L_{i-1}^k T_{i-1}E_{i}\\
\end{aligned} $$
Note that $ E_iE_{i-1} T_i L_{i-1}^k T_{i-1}E_{i}=E_iE_{i-1}
L_{i-1}^k E_{i-1}E_{i}=\omega_{i-1}^{(k)} E_i$ and
$E_iE_{i-1}E_iL_{i-1}^k T_{i-1}E_{i} =r L_{i-1}^k E_i$. By
(\ref{center}), together with (\ref{f1}), we have $E_iL_i^k
E_i=\omega_i^{(k)} E_i$, where $\omega_i^{(k)}\in R[L_1^\pm, L_2^\pm
, \cdots, L_{i-1}^\pm]$. We close the proof by showing that
$\omega_i^{(k)}\in Z(\cba{i-1})$. Note that any element $h\in
\cba{i-1} $ commutes with $E_i$ and $L_i$. We have $ h E_{i}L_i^k
E_i=E_{i}L_i^k E_i h$, which implies $E_i h\omega_i^{(k)} =E_{i}
\omega_i^{(k)} h$. If we identify the monomials of $\cba{i+1}$ with
Kauffman's tangles, we have $h E_{i}=0$ for $h\in \cba{i-1}$ if and
only if $h=0$. Thus $h\omega_i^{(k)} =\omega_i^{(k)} h$ for all
$h\in \cba{i-1}$.
\end{proof}

In the remainder of this section, we are going to construct the
 Jucys-Murphy basis of $\cba{n}$. We start by recalling
 some combinatorics.

Recall that a \textsf{partition} of $n$ is a weakly decreasing
sequence  of non--negative integers
$\lambda=(\lambda_1,\lambda_2,\dots)$ such that
$|\lambda|:=\lambda_1+\lambda_2+\cdots=n$. In this case, we write
$\lambda\vdash n$. The set $\Lambda^+(n)$, which consists of all
partitions of $n$, is a poset with dominance order $\trianglelefteq$
as the partial order on it. Given $\lambda, \mu\in \Lambda^+(n)$,
$\lambda\trianglelefteq \mu$ if  $ \sum_{j=1}^i \lambda_j\le
\sum_{j=1}^i \mu_j$ for all possible $i$. Write
$\lambda\vartriangleleft \mu$ if $\lambda\trianglelefteq \mu$ and
$\lambda\ne \mu$.

Suppose that  $\lambda$ and $\mu$ are two partitions. We say that
$\mu$ is obtained from $\lambda$ by \textsf{adding} a box if there
exists an $i$  such that $\mu_i=\lambda_i+1$ and $\mu_j=\lambda_j$
for $j\neq i$. In this situation we will also say that $\lambda$ is
obtained from $\mu$ by \textsf{removing} a box and we write
$\lambda\rightarrow\mu$ and $\mu\setminus\lambda=(i,\lambda_i+1)$.
We will also say that the pair  $(i,\lambda_i+1)$ is an
\textsf{addable} node of $\lambda$ and a \textsf{removable} node of
$\mu$. Note that $|\mu|=|\lambda|+1$.

The Young diagram $Y(\lambda)$ for a partition $\lambda=(\lambda_1,
\lambda_2, \cdots)$ is a collection of boxes arranged in
left-justified rows with $\lambda_i$ boxes in the $i$-th row of
$Y(\lambda)$. A $\lambda$-tableau $\ss$ is obtained by inserting $i,
1\le i\le n$ into $Y(\lambda)$ without repetition. The symmetric
group $\mathfrak S_n$ acts on $\ss$ by permuting its entries. Let
$\ts^\lambda$ be the $\lambda$-tableau obtained from the Young
diagram $Y(\lambda)$ by adding $1, 2, \cdots, n$ from left to right
along each row and from top to bottom along each column. If
$\ts^\lambda w=\ss$, write $w=d(\ss)$. Note that $d(\ss)$ is
uniquely determined by $\ss$.

A $\lambda$-tableau $\ss$ is standard if the entries in $\ss$ are
increasing both from left to right in each row and from top to the
bottom in each column. Let $\Std_n(\lambda)$ be the set of all
standard $\lambda$-tableaux.

Given an  $\ss\in \Std_n(\lambda)$, let $\ss\!\!\downarrow_i$ be
obtained from $\ss$ by removing all the entries $j$ in $\ss$ with
$j>i$. Let $\s_i$ be the partition of $i$ such that
$\ss\!\!\downarrow_i$ is an $\s_i$-tableau. Then  $\s=(\s_0, \s_1,
\cdots, \s_n)$ is a sequence of partitions such that
$\s_i\rightarrow \s_{i+1}$. Conversely, if we insert $i$ into the
box $\s_i\setminus \s_{i-1}$, then we obtain  an $\ss\in
\Std_n(\lambda)$. Thus, there is a bijection between
$\Std_n(\lambda)$ and the set of all $(\s_0,\s_1, \cdots, \s_n)$
such that  $\s_i\rightarrow \s_{i+1}, 0\le i\le n-1$ and $\s_0=0$,
and $\s_n=\lambda$.

Recall that $\mathfrak S_n$ is generated by $s_i, 1\le i\le n-1$
subject to the relations (1) $s_i^2=1$, $1\le i\le n-1$ (2)
$s_is_j=s_js_i$ if $|i-j|>1$ (3) $s_is_{i+1}s_i=s_{i+1}s_is_{i+1},
1\le i\le n-2$.
 Assume that  $0\le f\le
\lfloor\frac n2 \rfloor$. Let $\mathfrak S_{n-2f}$ be the subgroup
of $\mathfrak S_n$ generated by $s_j$, $2f+1\le j\le n-1$. Following
\cite{Enyang}, let $\mathfrak B_f$ be the subgroup of $\mathfrak
S_n$ generated by $\tilde s_{i}, \tilde s_0$, where $\tilde
s_i=s_{2i} s_{2i-1}s_{2i+1}s_{2i} $, $1\le i\le f-1$ and $\tilde
s_0=s_1$. Enyang~\cite{Enyang} proved that
 $\mathcal D_{f, n}$  is  a  complete set of right coset
representatives of $\mathfrak B_f\times \mathfrak S_{n-2f}$ in
$\mathfrak S_n$, where
$$\mathcal D_{f, n}=\Set[85]w\in \Sym_n| $(2i+1)w<(2j+1)w$,
$(2i+1)w<(2i+2)w$, $ 0\le
i<j<f$, and $ (k)w<(k+1) w$,  $2f< k< n$|.
$$
For $\lambda\vdash n-2f$, let $\mathfrak S_{\lambda}$ be the Young
subgroup of $\mathfrak S_{n-2f}$ generated by $s_j$, $2f+1\le j\le
n-1$ and $j\neq 2f+\sum_{k=1}^i \lambda_k$ for all possible $i$. A
standard $\lambda$-tableau $\hat \ss$ is obtained by using $2f+i,
1\le i\le n-2f$ instead of $i$ in the usual standard
$\lambda$-tableau $\ss$. Define $d(\hat \ss)\in \mathfrak S_{n-2f}$
by declaring that $\hat \ss=\hat{\ts}^\lambda d(\hat \ss)$. By abuse
of notation, we denote by $\Std_n(\lambda)$ the set of all standard
$\lambda$-tableaux $\hat\ss$.

It has been proved in ~\cite{Xi} that $\cba{n}$  is a cellular
algebra over a commutative ring. In what follows, we recall Enyang's
cellular basis for $\cba{n}$.

Let $\Lambda_n=\set{(f, \lambda)\mid \lambda\vdash n-2f, 0\le f\le
\lfloor\frac n 2\rfloor}$. Given  $(k, \lambda),  (f, \mu)\in
\Lambda_n$, define  $(k, \lambda)\trianglelefteq (f, \mu)$ if either
$k<f$ or $k=f$ and $\lambda\trianglelefteq\mu$. Write $(k,
\lambda)\lhd(f, \mu)$, if  $(k, \lambda)\unlhd (f, \mu)$ and $(k,
\lambda)\ne (f, \mu)$.

For any $w\in \mathfrak S_n$, write $T_w=T_{i_1}T_{i_2}\cdots
T_{i_k}$ if $s_{i_1}\cdots s_{i_k}$ is a reduced expression of $w$.
It is well-known that $T_w$ is independent of a  reduced expression
of $w$. Let $I(f, \lambda)=\Std_n(\lambda)\times \mathcal D_{f, n}$
and define
\begin{equation}
 C_{(\ss, u) (\ts, v)}^{(f, \lambda)} =T_u^\ast T_{d(\ss)}^{\ast}
  \mathfrak M_\lambda
T_{d(\ts)}T_v,\quad (\ss, u), (\ts, v)\in I(f, \lambda)
\end{equation}
where $\mathfrak M_\lambda=E^f X_{\lambda}$, $E^f=E_1 E_3\cdots
E_{2f-1}$, $X_\lambda=\sum_{w\in \mathfrak S_\lambda} q^{l(w)}T_ w$,
and $l(w)$,  the length of $w\in \mathfrak S_n$.

\begin{Theorem}\cite{Enyang} \label{cell} Let $\cba{n}$ be the Birman-Wenzl
algebra over  $R$. Let $\ast : \cba{n}\rightarrow \cba{n}$ be the
$R$-linear anti-involution in Lemma~\ref{anti}. Then
\begin{enumerate}
\item $\mathscr C_n =\left\{ C_{(\ss, u) (\ts, v)}^{(f, \lambda)} \mid (\ss, u), (\ts, v)\in
I(f, \lambda), \lambda\vdash n-2f,  0\le f\le \lfloor\frac
n2\rfloor\right\}$ is a free $R$--basis of $\cba{n}$.
\item $\ast (C_{(\ss, u) (\ts, v)}^{(f, \lambda)})=C_{(\ts, v) (\ss,
u)}^{(f, \lambda)}$.
\item For any $h\in \cba{n}$,
$$C_{(\ss, u) (\ts, v)}^{(f, \lambda)} h \equiv \sum_{(\us, w)\in I(f, \lambda)} a_{\us,
w} C_{(\ss, u) (\us, w)}^{(f, \lambda)} \mod
\cba{n}^{\vartriangleright (f, \lambda)}$$
 where
$\cba{n}^{\vartriangleright (f, \lambda)}$ is the free $R$-submodule
generated by $ C_{(\tilde \s, \tilde u) (\tilde \ts, \tilde v)}^{(k,
\mu)}$ with $(k, \mu) \vartriangleright (f, \lambda)$ and $(\tilde
\ss, \tilde u), (\tilde \ts, \tilde v)\in I(k, \mu)$.
 Moreover, each coefficient $a_{\us, w}$ is independent of $(\ss,
u)$.
\end{enumerate}
\end{Theorem}

Theorem~\ref{cell} shows that $\mathscr C_n $ is a cellular basis of
$\cba{n}$ in the sense of \cite{GL}. In this paper, we will only
consider right modules.

By  general theory about cellular algebras in \cite{GL}, we know
that,  for each $(f, \lambda)\in \Lambda_n$,  there is a cell module
$\Delta(f, \lambda)$ of $\cba{n}$, spanned by
$$\set{ \M_\lambda T_{d(\ts)}T_{v} \mod \cba{n}^{\vartriangleright
(f, \lambda)}\mid (\ts, v)\in I(f, \lambda)}.$$

We need Enyang's basis for $\Delta(f, \lambda)$ which is indexed
by up-down tableaux.

Given a $(f, \lambda)\in \Lambda_n$. An \textsf{$n$--updown
$\lambda$--tableau}, or more simply an updown $\lambda$--tableau, is
a sequence $\t=(\t_0, \t_1,\t_2,\dots,\t_n)$ of partitions such that
$\t_n=\lambda$, $\t_0=\varnothing$,  and either $\t_{i-1}\rightarrow
\t_{i}$ or $\t_{i}\rightarrow \t_{i-1}$  for $i=1,\dots,n$. Let
$\UPD_n(\lambda)$ be the set of updown $\lambda$--tableaux of $n$.

In what follows, we define  $T_{i, j}=T_{i}T_{i+1}\cdots T_{j-1}$
(resp. $T_{i-1}T_{i-2}\cdots T_j$) if $j>i$ (resp. if $j<i$). If
$i=j$, we set $T_{i, j}=1$.

\begin{Defn}\label{def-m}(cf. \cite{Enyang}) Given  $\t\in \UPD_n(\lambda)$ with
$\lambda \in \Lambda^+(n-2f)$, $0\le f\le \lfloor\frac n2\rfloor$,
define the non-negative integer $f_j$, $1\le j\le n$ and $0\le
f_j\le \lfloor j/2\rfloor$ by declaring that $\t_j\vdash j-2f_j$.
Let $\mu^{(j)}=\t_j$. Define $\M_\t=\M_{\t_n}$ inductively by
declaring that
\begin{itemize}\item[(1)]
$\M_{\t_1}=1$, \item[(2)]  $ \M_{\t_i}=\sum_{j=a_{k-1}+1}^{a_k}
q^{a_k-j} T_{j, i} \M_{\t_{i-1}}$ if $\t_i=\t_{i-1}\cup p$ with
$p=(k, \mu^{(i)}_k)$, and $a_l=2f_i+\sum_{j=1}^\ell \mu^{(i)}_j$
\item[(3)] $\M_{\t_i}= E_{2f_i-1} T_{i, 2f_i}^{-1}  T_{b_k,
2f_i-1}^{-1} \M_{\t_{i-1}} $ if $\t_{i-1}=\t_{i}\cup p$ with
$p=(k, \mu_k^{(i-1)})$,  and $b_k=2(f_i-1)+\sum_{j=1}^k
\mu^{(i-1)}_j$.\end{itemize}
\end{Defn}

It follows from  the definition that $\M_{\t}=\M_\lambda b_{\t}$ for
some $b_{\t}\in \cba{n}$. The following recursive formula describe
explicitly the element $b_{\t}$. Note that $b_\t=b_{\t_n}$ and
$\t_{n-1}=\mu$.
\begin{equation}\label{des-b}
b_{\t_{n}}=\begin{cases} T_{a_k, n}b_{\t_{n-1}}, & \text{ if
$\t_{n}= \t_{n-1}\cup\{ (k, \lambda_k)\}$}\\
T_{ n, 2f}^{-1}\sum_{j=b_{k-1}+1}^{b_k} q^{b_k-j} T_{j, 2f-1}^{-1}
b_{\t_{n-1}}, &
\text{ if $\t_{n-1}= \t_{n}\cup\{ (k, \mu_k)\}$}.\\
\end{cases}
\end{equation}

Suppose $\lambda\in \Lambda^+(n-2f)$ with $s$ removable nodes $p_1,
p_2, \cdots, p_s$ and $m-s$ addable nodes $p_{s+1}, p_{s+2}, \dots,
p_m$.
\begin{itemize}\item
Let $\mu^{(i)}\in \Lambda^+(n-2f-1)$  be obtained from $\lambda$ by
removing the  box $p_i$  for $1\le i\le s$.
\item Let  $\mu^{(j)}\in
\Lambda^+(n-2f+1)$  be obtained from $\lambda$ by adding the box
$p_{j}$ for  $s+1\le j\le m$.
\end{itemize}

 We identify  $\mu^{(i)}$ with $(k_i, \mu^{(i)})\in
\Lambda_{n-1}$  for $1\le i\le m$. So, $\mu^{(i)}\rhd\mu^{(j)}$ for
each $i, j$ with $1\le i\le s$ and $s+1\le j\le m$, and $k_i=f$ if
$1\le i\le s$ and $f-1$ otherwise. We arrange  $(k_i, \mu^{(i)})$'s
such that $(k_1,\mu^{(1)})\vartriangleright (k_2, \mu^{(2)})\cdots
\vartriangleright (k_m, \mu^{(m)})$.

 Define
$$
\begin{aligned} N^{\trianglerighteq \mu^{(i)}}& =R\text{--span} \{\M_\t\pmod
{\cba{n}^{\rhd (f, \lambda)}}\mid \t\in
\UPD_n(\lambda), \t_{n-1}\unrhd \mu^{(i)}\},\\
N^{\rhd \mu^{(i)}}& =R\text{--span} \{\M_\t\pmod {\cba{n}^{\rhd (f,
\lambda)}}\mid \t\in
\UPD_n(\lambda), \t_{n-1}\rhd \mu^{(i)}\}.\\
\end{aligned}
$$
In order to simplify the notation, we use $\bar \M_\t$ instead of
$\M_\t \pmod {\cba{n}^{\vartriangleright (f, \lambda)}}$ later on.
The following result is due to Enyang.
\begin{Theorem}\cite{Enyang}\label{ud} Let $\cba{n}$ be the Birman-Wenzl algebra
over $R$. Assume that $(f, \lambda)\in \Lambda_n$.
 \begin{enumerate}
 \item
$\{\bar \M_\t\mid \t\in \UPD_n(\lambda)\}$ is an $R$-basis of
$\Delta(f, \lambda)$.
\item Both $N^{\trianglerighteq \mu^{(i)}}$ and  $N^{\rhd
\mu^{(i)}}$ are $\cba{n-1}$-submodules of $\Delta(f, \lambda)$.
\item The $R$-linear map $\phi: N^{\unrhd \mu^{(i)}}/N^{\rhd
\mu^{(i)}} \rightarrow  \Delta(k_i, \mu^{(i)})$ sending
$\M_{\t}\pmod { N^{\rhd \mu^{(i)}}} $ to $ \M_{\t_{n-1}}\pmod{
\cba{n-1}^{\rhd (k_i, \mu^{(i)})}}$ is an isomorphism of
$\cba{n-1}$-modules.
\end{enumerate}\end{Theorem}

\begin{Defn} \label{udbasis} Given  $\s, \t\in \UPD_n(\lambda)$, define $\M_{\s,
\t}=b_{\s}^\ast \M_{\lambda}  b_{\t}$ where $\ast:
\cba{n}\rightarrow \cba{n} $ is the $R$-linear anti-involution on
$\cba{n}$ defined in Lemma~\ref{anti}.
\end{Defn}

Standard arguments prove the following result (cf.
\cite[Theorem~2.7]{RS:gram}).

\begin{Cor}\label{Mur}  Suppose that $\cba{n}$ is  the Birman-Wenzl algebra
over  $R$. Then
\begin{enumerate}
\item $\mathscr M_n =\{ \M_{\s \t}\mid \s,  \t\in \UPD_n(\lambda), \lambda
\vdash n-2f, 0\le f\le \floor{ n2}\} $ is a free $R$--basis of
$\cba{n}$.
\item $ \M_{\s \t}^\ast =\M_{\t \s}$ for all $\s, \t\in \UPD_n(\lambda)$
and all $(f, \lambda)\in \Lambda_n$.
\item  Let $\widetilde {\cba{n}}^{\vartriangleright (f, \lambda)}$ be the free
$R$-submodule of $\cba{n}$ generated by $ \M_{\tilde \s \tilde
\t}$ with $\tilde \s, \tilde \t\in \UPD_n(\mu)$ and $(\frac
{n-|\mu|} 2, \mu) \vartriangleright (f, \lambda)$. Then
$\widetilde {\cba{n}}^{\vartriangleright (f,
\lambda)}=\cba{n}^{\vartriangleright (f, \lambda)}$.
\item For all $\s, \t\in \UPD_n(\lambda)$, and all $h\in \cba{n}$,
there exist scalars $a_\u\in R$ which are independent of $\s$,
such that
$$\M_{\s \t} h \equiv \sum_{\u} a_{\u} \M_{\s\u} \pmod{
 {\cba{n}}^{\vartriangleright (f, \lambda)}}.$$
\end{enumerate}
\end{Cor}

We call $\mathscr M_n$ the \textsf{Jucys-Murphy} basis of $\cba{n}$.
It is a cellular basis of $\cba{n}$ over $R$. In \cite{GL}, Graham
and Lehrer proved that there is a symmetric invariant bilinear form
$\langle\quad, \quad \rangle: \Delta(f, \lambda)\times  \Delta(f,
\lambda)\rightarrow R$  on each cell module. In our case, we use
$\mathscr M_n$ to define such a bilinear form on $\Delta(f,
\lambda)$. More explicitly, $\langle \bar \M_\s, \bar \M_\t \rangle
\in R $ is determined by
$$\M_{\tilde\s\s}\M_{\t\tilde\t}\equiv\langle \bar \M_\s,  \bar \M_\t   \rangle
\M_{\tilde\s\tilde\t}\pmod { \cba{n}^{\vartriangleright (f,
\lambda)}},\quad  \tilde\s, \tilde \t\in \UPD_n(\lambda).$$

By  Corollary~\ref{Mur}(d), the above symmetric invariant bilinear
form is independent of $\tilde\s, \tilde\t\in \UPD_n(\lambda)$. The
Gram matrix $G_{f, \lambda}$ with respect to the Jucys-Murphy basis
of $\Delta(f, \lambda)$ is the $k\times k$ matrix with $$k=\rank
\Delta(f, \lambda)=\frac{n!(2f-1)!!}{(2f)! \prod_{(i, j)\in \lambda}
h_{i, j}^\lambda}$$ where $h_{i,
j}^\lambda=\lambda_i+\lambda_j'-i-j+1$ is the hook length. The $(\s,
\t)$-th entry of $G_{f, \lambda}$ is $\langle \bar \M_\s,  \bar
\M_\t \rangle$.

One of the main purposes of this paper is to compute the Gram
determinant $\det G_{f, \lambda}$ associated to each cell module
$\Delta(f, \lambda)$.

Given $\s\in \UPD_n(\lambda)$, we identify $\s_i$ with $(f_i,
\mu^{(i)})$ if $\s_i=\mu^{(i)}\vdash i-2f_i$. Define the partial
order $\trianglelefteq $ on $\UPD_n(\lambda)$ by declaring that
$\s\trianglelefteq\t$ if $\s_i\trianglelefteq\t_i$ for all $1\le
i\le n$. Write $\s\triangleleft\t$ if $\s\trianglelefteq\t$ and
$\s\neq \t$. We remark that Enyang has used $\unlhd$ to state
Theorem~\ref{xproduct}.  We define the relation  $\succ$ on
$\UPD_n(\lambda)$ instead of his partial order $\unlhd$.

Suppose $\s\neq \t$.  We write $\s\succ\t$   if there is a positive
integer $ k\le n-1$ such that $\s_k\rhd \t_k$ and $\s_j=\t_j$ for
$k+1\le j\le n$. We will use $\s\overset{k} \succ\t$ to denote
$\s_j\rhd\t_j$ and $\s_\ell=\t_\ell$ for $j+1\le \ell\le n$ and
$j\ge k$.

For any $(f, \lambda)\in \Lambda_n$, define  $\t^\lambda\in
\UPD_n(\lambda)$ such that
\begin{itemize}
\item $\t^\lambda_{2i-1}=(1)$ and $\t^\lambda_{2i}=\varnothing$ for $1\le i\le f$,
\item   $\t^\lambda_{i}$ is obtained from $\hat
\ts^\lambda$ by removing the entries $j$ with $j>i$ under the
assumption  $2f+1\le i\le n$.
\end{itemize}
Then $\t^\lambda$ is maximal in $\UPD_n(\lambda)$ with respect to
$\succ$ and $\unrhd$.

For any  $\t\in \UPD_n(\lambda)$ with $(f, \lambda)\in \Lambda_n$,
define $c_\t(k)\in R$ by
$$c_\t(k)= \begin{cases} \label{content}
   rq^{2(j-i)},  &\text{if }\t_k=\t_{k-1}\cup (i,j),\\
     r^{-1}q^{2(i-j)}, &\text{if }\t_{k-1}=\t_{k}\cup(i,j).
\end{cases}$$
 If $p=(i, j)$ is an addable (resp. a
removable)  node of $\lambda$, define $c_\lambda(p)=  j-i$ (resp.
$- j+i $).

The following result plays a key role in the construction of an
orthogonal basis for $\cba{n}$.

\begin{Theorem}\cite{Enyang} \label{xproduct} Given $\s, \t\in \UPD_n(\lambda)$,
 with $(f, \lambda)\in \Lambda_n$,   $$\M_{\s\t} L_k\equiv c_\t(k)\M_{\s\t} +
\sum_{{\u\overset{k-1}\succ\t}} a_\u \M_{\s\u} \pmod
{\cba{n}^{\vartriangleright (f, \lambda)}} .$$\end{Theorem}

Enyang used $\u\rhd\t$ instead of $\u\succ\t$.  However, we could
not understand  the claim about $\breve{N}^\mu$ under
\cite[(7.3)]{Enyang}. If one uses $\succ$ instead of $\rhd$, then
everything in the proof of \cite[7.8]{Enyang} is available.

\section{orthogonal representations for $\cba{n}$} In this section, we assume that $F$ is a  field which
contains non-zero $q, r$ and $(q-q^{-1})^{-1}$ such that $o(q^2)>n$
and $|c|> 2n-3$ whenever $ r^2 q^{2c}=1$ for some $c\in \mathbb Z$.
 The main purpose of this section is to construct an orthogonal
basis of $\cba{n}$ over $F$.

Suppose $1\le k\le n$ and $(f, \lambda)\in \Lambda_n$. Define an
equivalence relation $\simk$ on $\UPD_n(\lambda)$ by declaring
that $\t\simk \s$ if $\t_j=\s_j$ whenever $1\le j\le n$ and $j\neq
k$, for $\s,\t\in\UPD_n(\lambda)$. The following result is
well-known. See e.g. \cite{RS:gram}.

\begin{Lemma} Suppose $s\in \UPD_n(\lambda)$ with
$\s_{k-1}=\s_{k+1}$. Then there is a bijection between the set of
all addable and removable nodes of $\s_{k+1}$ and the set
$\set{\t\in \UPD_n(\lambda)\mid \t\simk \s}$.
\end{Lemma}

Suppose $\lambda$ and $\mu$ are partitions. We write
$\lambda\ominus\mu=\alpha$ if either $\lambda\supset \mu$ and
$\lambda\setminus \mu=\alpha$ or $\lambda\subset \mu$ and
$\mu\setminus \lambda=\alpha$. The following lemma
 can be proved by
  arguments similar to those in \cite{RS:gram}.

\begin{Lemma} \label{cont}
Assume  that $\s, \t\in \UPD_n(\lambda)$ with $(f, \lambda)\in
\Lambda_n$.
\begin{enumerate} \item $\s=\t$ if and only if $c_\s(k)=c_\t(k)$
for $1\le k\le n$.
\item Suppose $\t_{k-1}\neq \t_{k+1}$. Then  $c_\t(k)\neq  c_\t(k+1)$.
\item If $\t_{k-1}=\t_{k+1}$, then $c_\t(k)\neq c_\s(k)^{\pm}$
whenever $\s\simk\t$ and $\s\neq \t$.
\item $c_\t(k)\not\in \{-q, q^{-1}\}$ for all $\t\in \UPD_n(\lambda)$ with
$\t_{k-1}=\t_{k+1}$.
\end{enumerate}
\end{Lemma}

\begin{Defn}\label{def-r} Suppose $\t\in \UPD_n(\lambda)$ for some $(f, \lambda)\in \Lambda_n$.
Following \cite{M:gendeg}, we define
 \begin{enumerate}\item
$ \mathscr R(k)=\set{ c_\s(k)\mid \s\in \UPD_n(\lambda), (f,
\lambda)\in \Lambda_n}, \quad 1\le k\le n$,
\item $F_\t =\prod_{k=1}^n F_{\t, k} $ where $F_{\t, k}=\prod_{\substack{r\in \mathscr R(k)\\
c_\t(k)\neq r}} \frac {L_k-r} {c_\t(k)-r}, \text{ and } f_\t= \bar
\M_\t F_\t, \quad \t\in \UPD_n(\lambda)$,
\item
$f_{\s\t}=F_\s \M_{\s\t} F_\t, \quad \s, \t\in
\UPD_n(\lambda).$\end{enumerate}
\end{Defn}

Standard arguments prove Lemma~\ref{fxproduct} and Lemma~\ref{epro}.
Mathas has proved similar results for a general class of cellular
algebras in \cite{M:semi}. Although he has used a partial order
which is similar to $\unlhd$, his arguments can be used to verify
the following results. See also \cite{M:ULect} for the Hecke algebra
$\H_n$ of type $A$.

\begin{Lemma}\label{fxproduct}
 Suppose that  $\t\in \UPD_n(\lambda)$ with  $(f, \lambda)\in \Lambda_n$.
\begin{enumerate}\item  $f_\t=\bar\M_\t + \sum_{\s\succ \t} a_\s
\bar \M_\s$.
 \item $\bar \M_\t=f_\t + \sum_{\s\succ \t} b_\s
f_\s$.
\item $f_\t L_k=c_\t(k) f_\t$, for any $k$, $1\le k\le n$.
\item $f_\t  F_\s =\delta_{\s\t} f_\t$ for all $\s\in
\UPD_n(\mu)$ with $(\frac {n-|\mu|}{2}, \mu)\in \Lambda_n$.
\item $\{f_\t\mid \t\in \UPD_n(\lambda)\}$ is a basis of
$\Delta(f, \lambda)$.
\item The Gram determinants associated to $\Delta(f, \lambda)$
defined by $\{f_\t\mid \t\in \UPD_n(\lambda)\}$ and $\{\bar
\M_\t\mid \t\in \UPD_n(\lambda)\}$ are the same.
\end{enumerate}
\end{Lemma}
 Let  $f_\t T_k=\sum_{\s\in \UPD_n(\lambda)} s_{\t\s}(k)
f_\s$ and $f_\t E_k=\sum_{\s\in \UPD_n(\lambda)} E_{\t\s}(k) f_\s$.
\begin{Lemma}\label{epro}
Suppose $\t\in \UPD_n(\lambda)$ and $1\le k\le n-1$.
\begin{enumerate}
\item $\s\simk \t$ if either  $s_{\t\s}(k)\neq 0$ or
$E_{\t\s}(k)\neq 0$. \item  $f_\t E_k=0$ if $\t_{k-1} \neq
\t_{k+1}$.
\item If $\t_k\ominus \t_{k-1}$ and $\t_{k+1} \ominus \t_k$ are
neither in the same row nor in the same column, then there is a
unique up-down tableau in $\UPD_n(\lambda)$, denoted by $\t s_k $,
such that  $\t s_k\simk \t$ and $c_{\t}(k)=c_{\t s_k}(k+1)$ and
$c_{\t}(k+1)=c_{\t
 s_k}(k)$.\item  If $\t_k\ominus \t_{k-1}$ and $\t_{k+1} \ominus
\t_k$ are either in the same row or in the same column, then there
is no $\s\in \UPD_n(\lambda)$ such that $\s\simk \t$ and
$c_{\t}(k)=c_{\s}(k+1)$ and $c_{\t}(k+1)=c_{\s}(k)$.
\end{enumerate}\end{Lemma}

\begin{Lemma}\label{mts}
Suppose that $\t\in\UPD_n(\lambda)$ with $\t_{i-2}\neq \t_i$, $\t
s_{i-1}\in \UPD_n(\lambda)$  and $\t s_{i-1}\lhd\t$. We have
\begin{itemize}
\item[a)] If $\t_{i-2}\subset\t_{i-1}\subset\t_i$, then
$\bar \M_{\t}T_{i-1}=\bar \M_{\t s_{i-1}}$.

\item[b)] If $\t_{i-2}\supset\t_{i-1}\subset\t_i$ such that
$\ell>k$ where $\t_{i-2}\setminus\t_{i-1}=(k,\nu_k)$,
$\t_{i}\setminus\t_{i-1}=(\ell,\mu_{\ell})$, $\t_{i-2}=\nu$ and
$\t_i=\mu$, then $\bar \M_{\t}T_{i-1}^{-1}=\bar \M_{\t s_{i-1}}$.
\end{itemize}
\end{Lemma}
\begin{proof} The proof of the result is essentially identical to the proof of the corresponding
result in the proof of \cite[3.14]{RS:gram}. One can check it by
Definitions~\ref{bmw-def} and \ref{def-m}. We leave the details to
the reader.
\end{proof}

\begin{Lemma}\label{fss}
 Suppose $\t\in \UPD_n(\lambda)$ with $\t_{i-1} \neq \t_{i+1}$ and
 $\t s_i\in \UPD_n(\lambda)$. Then $f_\t T_i=s_{\t\t}(i) f_\t + s_{\t,\t s_i} (i) f_{\t s_i} $,
where
 \begin{itemize} \item $s_{\t\t}(i) =\frac {\omega c_\t
(i+1)}{ c_\t(i+1)-c_\t(i)}$, \item $ s_{\t,\t s_i}(i) = 1-\frac
{c_\t(i)}{c_\t (i+1)}s_{\t\t}^2(i)$ if
 $\t s_i \rhd \t$ and $ s_{\t,\t s_i}(i) = 1$ if $\t s_i \lhd \t$ and  one of
 the following conditions holds,
\begin{enumerate} \item  $\t_{i-1}\subset\t_{i}\subset\t_{i+1}$,
\item $\t_{i-1}\supset\t_{i}\subset\t_{i+1}$ such that $\ell>k$ where
$\t_{i-1}\setminus\t_{i}=(k,\nu_k)$,
$\t_{i+1}\setminus\t_{i}=(\ell,\mu_{\ell})$, $\t_{i-1}=\nu$ and
$\t_{i+1}=\mu$.\end{enumerate} \end{itemize}
\end{Lemma}
\begin{proof}
Write $f_\t T_i=\sum_{ \s\overset {i}\sim\t} s_{\t\s}(i) f_\s$. By
Lemma~\ref{epro}, $s_{\t\s}(i)\neq 0$ implies $\s\in \{\t, \t
s_i\}$, and $f_\t E_i =0$. On the other hand, by
Lemma~\ref{wenzl-rel}(c),
$$ T_i L_{i+1} =L_i T_i+\omega L_{i+1}-\omega r^{-1} E_i L_i T_i.$$
So, $f_\t T_i L_{i+1}=c_\t(i) f_\t T_i +\omega c_\t(i+1) f_\t.$
Comparing the coefficient of $f_\t$ in   $f_\t T_i L_{i+1}$ yields
the formula on $s_{\t\t}(i)$ as required.

We compute $ s_{\t,\t s_i}(i)$ under the assumptions as follows. By
Lemma~\ref{mts}, $\bar \M_\t T_{i}=\bar \M_{\t s_{i}}$ if
$\t_{i-1}\subset\t_{i}\subset\t_{i+1}$ and $\t s_i\lhd \t$.
 By
Lemma~\ref{fxproduct}(a)--(b),
$$ f_\t T_{i} =(\bar \M_\t
+\sum_{\u\succ \t} a_\u f_\u) T_{i} =\bar \M_{\t s_{i}} +
\sum_{\u\succ \t} a_\u f_\u T_{i} .$$
 If $f_{\t s_{i}}$
appears in the expression of $f_\u T_{i}$ with non-zero coefficient,
then $\u\overset {i}  \sim \t s_{i}$. Therefore, $\u\in \{ \t, \t
s_{i}\}$ which contradicts $\u\succ \t\rhd  \t s_{i}$. By
Lemma~\ref{fxproduct}(a), the coefficient of $f_{\t s_i}$ in $f_\t
T_{i}$ is $1$.

Under the assumption given in (b), we have $\bar
\M_{\t}T_{i}^{-1}=\bar \M_{\t s_{i}}$. Thus, $$\begin{aligned}
\bar\M_{\t} T_{i} & =\bar \M_{\t s_{i}} T_{i}^2=\bar \M_{\t s_i}
(1+\omega(T_{i}-r^{-1}
E_{i}))\\
&= \bar \M_{\t s_{i}}+\omega\bar \M_{\t} -\omega r^{-1}\bar \M_{\t
s_{i}} E_{i}\\
\end{aligned}
$$
By Lemma~\ref{fxproduct}, $\bar \M_{\t s_i}=f_{\t s_{i}}
+\sum_{\u\succ \t s_{i}} a_\u f_\u$ for some $a_\u\in F$. Since $\t
s_{i}\overset i\sim \t$, $(\t s_{i})_{i-1}\neq (\t s_{i})_{i+1}$. By
Lemma~\ref{epro}(b), $f_{\t s_{i}} E_{i}=0$. If $f_{\t s_{i}}$
appears in the expression of $f_\u E_{i}$ with non-zero coefficient,
then $\u\overset i \sim  \t s_{i}$, forcing $\u_{i-1}\neq \u_{i+1}$.
Thus, $f_\u E_{i}=0$, a contradiction. Finally, by
Lemma~\ref{fxproduct}(a), the coefficient of $f_{\t s_{i}}$ in
$\bar\M_{\t s_i} $ is $1$, forcing $s_{\t, \t s_{i}}(i)=1$.
\end{proof}

Note that the bilinear form $\langle\ , \ \rangle: \Delta(f,
\lambda)\times \Delta(f, \lambda)\rightarrow F$ is associative. We
have $\langle f_\t T_k, f_\t T_k  \rangle =\langle f_\t T_k^2, f_\t
\rangle$.  The proofs of Corollary~\ref{up1} and Lemma~\ref{samerow}
are essentially identical to \cite[4.3, 3.15]{RS:gram}. We leave the
details to the reader.
\begin{Cor}\label{up1}
Suppose $\t\in \UPD_n(\lambda)$ with $(f, \lambda)\in \Lambda_n$
and $\t_{k-1}\neq \t_{k+1}$. If $\t s_k\in \UPD_n(\lambda)$ and
$\t s_k \vartriangleleft \t$, then
$$
\langle f_{\t s_k}, f_{\t s_k} \rangle=\left(1-\frac{\omega^2
c_\t(k) c_\t(k+1)}{(c_\t(k+1)- c_\t(k))^2} \right) \langle f_\t,
f_\t \rangle.
$$
\end{Cor}

\begin{Lemma}\label{samerow}
 Suppose $(f, \lambda)\in \Lambda_n$. If $\t\in \UPD_n(\lambda)$
 with $\t_{k-1}\neq \t_{k+1}$, then
 \begin{enumerate}\item  $f_\t T_k=q f_\t$ if $\t_k\ominus \t_{k-1}$ and
 $\t_k \ominus \t_{k+1}$ are in the same row,
\item  $f_\t T_k=-q^{-1} f_\t$ if $\t_k\ominus \t_{k-1}$ and
 $\t_k \ominus \t_{k+1}$ are in the same column.
\end{enumerate}
\end{Lemma}

 In the following, we assume that $F=
\mathbb C(r^\pm, q^\pm, \omega^{-1})$, where $r, q$ are
indeterminates  and $\omega=q-q^{-1}$.

\begin{Lemma}\label{ek} Suppose that  $\t\in \UPD_n(\lambda)$ and
$\t_{k-1}=\t_{k+1}$. Then
\begin{enumerate}
\item $
    E_{\t\t}(k)=r c_\t(k)^{-1} (1+\omega^{-1} (c_\t(k)-c_\t(k)^{-1})    \prod_{\substack{\s\simk\t\\ \s\neq \t}}
    \frac{c_\t(k)-c_\s(k)^{-1}}{c_\t(k)-c_\s(k)}\neq 0$.

\item $E_{\t\s}(k)E_{\u\u}(k)=E_{\t\u}(k)E_{\u\s}(k)$ for any
$\s, \t, \u\in \UPD_n(\lambda)$ with $\t\simk\s\simk \u$.
\end{enumerate}
\end{Lemma}
\begin{proof}
First, we prove $E_{\t\t}(k)\neq 0$. By assumption and  \cite[
5.6]{W2},  $\Delta(f, \lambda)$ is  irreducible since  $\cba{n}$ is
semisimple. In \cite[6.17]{LR}, Leduc and Ram proved   that the
seminormal representations $S^{f, \lambda}$\footnote{In \cite{LR},
$S^{f, \lambda}$ is denoted by $\mathscr Z^\lambda$} over $\mathbb
C$  with special parameters  for all $(f, \lambda)\in \Lambda_n$
consist of the complete set of pair-wise non-isomorphic irreducible
modules when $\cba{n}$ is semisimple. By the fundamental theorem of
algebra, one can get the same results over $\mathbb C(q^\pm, r^\pm,
\omega^{-1})$ where $r, q$ are indeterminates and $\omega=q-q^{-1}$.

Thus, $\Delta(f, \lambda)\cong S^{\ell, \mu}$ for some $(\ell,
\mu)\in \Lambda_n$. If we denote by $\phi$ the corresponding
isomorphism between $\Delta(f, \lambda)$ and $ S^{\ell, \mu}$, then
$\phi(f_\t)\in S^{\ell, \mu}$. In \cite{LR}, Leduc and Ram
constructed a basis for $S^{\ell, \mu}$, say $v_\s$, $\s\in
\UPD_n(\mu)$ such that $v_\s L_k=c_\s(k) v_\s$. Note that $f_\t\in
\Delta(f, \lambda) $ is a common eigenvector   of $L_k$, $1\le k\le
n$.  By Lemma~\ref{cont}(a), $(\ell, \mu)=(f, \lambda)$  and
$\phi(f_\t)$ is equal to $v_\t$ up to a scalar since the common
eigenspace on which $L_k, 1\le k\le n$ acts as $c_\s(k)$ is of one
dimension. Leduc and Ram \cite[5.9]{LR}\footnote{In \cite{LR},
$r=\epsilon q^a$ for some $a\in \mathbb Z$ and $\epsilon\in \{1,
-1\}$. Therefore, $\tilde E_{\t\t}(k)\neq 0$ if $r$ is an
indeterminant.}
 proved that $\tilde E_{\t\t}(k)\neq 0$ for any $\t\in
\UPD_n(\lambda)$, where $\tilde E_{\t\s}(k)$ is defined by $ v_\t
E_k=\sum_{\s\simk\t} \tilde E_{\t\s}(k) v_\s$. Since $v_\t$ is a
non-zero scalar of
 $\phi(f_\t)$,
 $E_{\t\t}(k)= \tilde E_{\t\t}(k)\neq 0$. We remark that
$E_{\t\t}(k)\in F$ since $f_\t$ is an $F$-basis element of
$\Delta(f, \lambda)$. In general, $E_{\t\s}(k)\neq \tilde
E_{\t\s}(k)$ if   $\s\neq \t$.

In \cite{BC}, Beliakova and Blanchet proved that the generating
function $W_k(y)=\sum_{a\ge 0} \omega_k^{(a)}/y^a$ satisfies the
following identity $$
\frac{W_{k+1}(y)+r^{-1}\omega^{-1}-\frac{y^2}{y^2-1}}
{W_{k}(y)+{r^{-1} \omega^{-1}}-\frac{y^2}{y^2-1}}=\frac{
y^{-1}-\frac{\omega^2 L_k^{-1}}{(y-L_k^{-1})^2}}{
y^{-1}-\frac{\omega^2 L_k}{(y-L_k)^2}}.
$$ Comparing the coefficients
of $f_\s$ on both sides of $f_\t E_k W_k(y)=f_\t E_k \frac{y}{y-L_k}
E_k$ yields
$$
{W_k(y, \s)}{y^{-1}} E_{\t\s}(k) =\sum_{\s\simk\t\simk \u}\frac{
E_{\t\u}(k) E_{\u\s} (k) } {y-c_\u(k)}.$$ Thus $ E_{\t\s}(k)\cdot
Res_{y=c_\u(k)} { W_k(y, \s)y^{-1}} =E_{\u\s}(k) E_{\t\u} (k)$.
Since $E_{\t\t}(k)\neq 0$,  $E_{\t\t}(k)=Res_{y=c_\t(k)}{ W_k(y,
\t)}{y^{-1}}$ by assuming that $\t=\s=\u$. By computation, we can
verify $$\underset{y=c_\t(k)} {\text{Res}}\frac{ W_k(y, \t)}{y}=r
c_\t(k)^{-1} (1+ \omega^{-1} (c_\t(k)-c_\t(k)^{-1}))
    \prod_{\substack{\s\simk\t\\ \s\neq \t}}
    \frac{c_\t(k)-c_\s(k)^{-1}}{c_\t(k)-c_\s(k)}.$$
This completes the proof of (a). If $\s\simk \u$, then
$c_\s(j)=c_\u(j)$ for $j\le k-1$. By Lemma~\ref{e-powerk},
$\omega_k^{(a)}\in F[L_1^\pm, L_2^\pm, \cdots, L_{k-1}^\pm]$. So,
${W_k(y, \s)}=W_k(y, \u)$. Thus
 $E_{\t\s}(k)E_{\u\u} (k) =E_{\u\s}(k) E_{\t\u} (k)$, proving (b).
\end{proof}

The following result is a special case of \cite[3.14]{M:semi} which
is about the construction of the primitive idempotents and central
primitive idempotents for a general class of  cellular algebras. In
our case, such idempotents can be computed explicitly via  a
recursive formula on $\langle f_\t, f_\t\rangle$, $\t\in
\UPD_n(\lambda)$ and $(f, \lambda)\in \Lambda_n$. This recovers the
main result in \cite{BC}.

\begin{Prop}\label{wed} Suppose that $\cba{n}$ is the
Birman-Wenzl algebra over  $\mathbb C(q^\pm, r^\pm, \omega^{-1})$.
\begin{enumerate}\item Suppose $\t\in \UPD_n(\lambda)$. Then
$\frac{1}{\langle f_\t, f_\t \rangle} f_{\t\t}$ is a  primitive
idempotent of semisimple  $\cba{n}$ with respect to the cell
module $\Delta(f, \lambda)$.
\item $\sum_{\t\in \UPD_n(\lambda)} \frac{1}{\langle f_\t, f_\t \rangle}
f_{\t\t}$ is a central primitive idempotent. Furthermore,
$$\sum_{(\frac{n-|\lambda|} 2, \lambda)\in \Lambda_n}
 \sum_{\t\in \UPD_n(\lambda)} \frac{1}{\langle f_\t,
f_\t \rangle} f_{\t\t}=1.$$
\end{enumerate}
\end{Prop}

\section{Gram determinants for $\cba{n}$}
In this section, we compute the Gram determinant for each cell
module of $\cba{n}$ over $F=\mathbb C(q^\pm, r^\pm, \omega^{-1})$,
where $r, q$ are indeterminates and $\omega=q-q^{-1}$.  Our result
for the Gram determinants holds true for $\cba{n}$ over $R:=\mathbb
Z[q^\pm, r^\pm, \omega^{-1}]$ since the Jucys-Murphy basis for
$\Delta(f, \lambda)$ is an $R$-basis. By base change, it holds over
an arbitrary field.

 Given $\t\in \UPD_n(\lambda)$ with $\t_{n-1}=\mu$, define $\hat
\t\in \UPD_{n-1}(\mu)$ such that $\hat\t_i=\t_i$, $1\le i\le n-1$,
and $\tilde \t\in \UPD_n(\lambda)$ with $\tilde\t_j=\t^\mu_j$ for
$1\le j\le n-1$ and $\tilde\t_n=\t_n=\lambda$.

Let  $[n]=1+q^2+\cdots +q^{2n-2}$ and $[n]!=[n][n-1]\cdots [2][1]$.
If $\lambda=(\lambda_1, \lambda_2, \cdots, \lambda_k)$, define
$[\lambda]!=\prod_{i=1}^k [\lambda_i]!$. Standard arguments prove
the following result (cf.\cite[4.2]{RS:gram}).
 We leave the details to the reader.
\begin{Prop} \label{recursive} Suppose $\t\in \UPD_n(\lambda)$ with $(f,
\lambda)\in \Lambda_n$. If $\t_{n-1}=\mu$ with $(\ell, \mu)\in
\Lambda_{n-1}$, then $ \langle  f_\t,  f_\t \rangle =\langle
f_{\hat\t} , f_{\hat\t} \rangle \frac {\langle f_{\tilde \t} ,
f_{\tilde\t} \rangle}{\delta^\ell [\mu]!}$.
\end{Prop}

\def\R{\mathscr R(\lambda)^{<p}}
\def\A{\mathscr A(\lambda)^{<p}}
\def\RR{\mathscr R(\mu)^{<p}}
\def\AA{\mathscr A(\mu)^{<p}}
\def\RB{\mathscr R(\lambda)^{\ge p}}
\def\AB{\mathscr A(\lambda)^{\ge p}}
\def\RRB{\mathscr R(\mu)^{\ge p}}
\def\AAB{\mathscr A(\mu)^{\ge p}}

For any $\lambda\vdash n-2f$, let $ \mathscr A(\lambda)$ (resp.
$\mathscr R(\lambda))$ be the set of all addable (resp. removable)
nodes of $\lambda$. Given a  $p=(k, \lambda_k)\in \mathscr
R(\lambda)$ (resp. $p=(k, \lambda_k+1)\in \mathscr A(\lambda)$,
define
\begin{enumerate}
\item  $\mathscr R(\lambda)^{<p} =\{(\ell, \lambda_\ell)\in \mathscr
R(\lambda) \mid \ell>k\}$,
\item $\mathscr A(\lambda)^{<p} =\{(\ell, \lambda_\ell+1)\in \mathscr
A(\lambda) \mid \ell>k\}$,
\item $\mathscr R(\lambda)^{\ge p} =\{(\ell, \lambda_\ell)\in \mathscr
R(\lambda) \mid \ell\le k\}$,
\item  $\mathscr A(\lambda)^{\ge p} =\{(\ell,
\lambda_\ell+1)\in \mathscr A(\lambda) \mid \ell\le  k\}$.
\end{enumerate}

\begin{Prop}\label{key0}
Suppose  $\t\in \UPD_n(\lambda)$ with $(f, \lambda)\in \Lambda_n$.
If $\hat \t=\t^\mu$ and   $\t_{n}= \t_{n-1}\cup \{p\}$ with $p=(k,
\lambda_k)$, then
\begin{equation}\label{key1}
\frac{\langle f_\t,  f_\t \rangle} {\delta^f [\mu]!}
=-q^{2\lambda_k} \frac{\prod_{r_1\in \A} [c_{\lambda}
(p)+c_{\lambda} (r_1)]}{\prod_{r_2\in \R} [c_{\lambda}
(p)-c_{\lambda} (r_2)]}.
\end{equation}
\end{Prop}

\begin{proof} By assumption, $\t=\t^\lambda s_{a, n}$ where
$a=2f+\sum_{j=1}^k \lambda_j$.  Note that  $\t\vartriangleleft\t
s_{n-1}\vartriangleleft\cdots\vartriangleleft \t s_{n, a}$. By
Proposition~\ref{recursive} and Corollary~\ref{up1},
\begin{equation}\label{up}\langle f_\t, f_\t \rangle=  \langle f_{\t^\lambda},
f_{\t^\lambda}\rangle \prod_{j=a+1}^n \left( 1-\omega^2
\frac{c_{\t^\lambda}(j)c_{\t^\lambda}(a)}
{(c_{\t^\lambda}(j)-c_{\t^\lambda}(a))^2}\right).\end{equation}

Since $f_{\t^\lambda}=\bar \M_{\t^\lambda}$, $\langle
f_{\t^\lambda}, f_{\t^\lambda}\rangle=\delta^f [\lambda]!$.
 Using the definitions of
$c_{\t^\lambda}(j)$ for $a\le j\le n$ to simplify (\ref{up}) yields
(\ref{key1}).
\end{proof}

\begin{Prop}\label{key21}  Suppose $\t \in \UPD_n(\lambda)$ with
$\lambda=(\lambda_1, \dots, \lambda_k)\vdash n-2f$. If $\t^{\mu}
=\hat \t$ and $\t_{n-1} = \t_n \cup p$ with $p=(k,\mu_k)$, then
\begin{equation}\label{key2} \frac{ \langle f_\t,\ f_\t\rangle }{\delta^{f-1}[\mu]!} =
[\mu_k] E_{\t\t}(n-1). \end{equation} \end{Prop}

\begin{proof}  Write $\widetilde
\M_\lambda=\M_{\t^\lambda}$ with $\t^\lambda\in
\UPD_{n-2}(\lambda)$. Let $a = 2(f-1)+\sum_{j=1}^{k-1} \mu_j+ 1$. By
Definition~\ref{def-m} and Lemma~\ref{wenzl-rel}(d)-(e),
 $$\begin{aligned}   f_\t E_{n-1}   \equiv &
E_{2f-1}T_{n,2f}^{-1} T_{n-1, 2f-1}^{-1} \widetilde
\M_{\lambda}\sum_{j=a}^{n-1} q^{n-1-j} T_{n-1,j} F_\t E_{n-1} \mod
\cba{n}^{\vartriangleright(f, \lambda)}\\
  \equiv & E_{2f-1}E_{2f}\cdots E_{n-1}\widetilde \M_\lambda \sum_{j=a}^{n-1}
q^{n-1-j} T_{n-1, j} F_{\t, n} F_{\t, n-1}\\
&\times  E_{n-1} \prod_{k=1}^{n-2} F_{\t, k} \mod
\cba{n}^{\vartriangleright(f, \lambda)}.\\
\end{aligned}$$
By Lemma~\ref{wenzl-rel}(b),  $E_{n-1}\widetilde\M_\lambda
=\widetilde\M_\lambda E_{n-1}$. Via Proposition~\ref{omega}, we
write
\begin{equation}\label{Phi}
E_{n-1}F_{\t,n-1}F_{\t, n}E_{n-1}=\Phi(L_1^\pm,\cdots,
L_{n-2}^\pm)E_{n-1}\end{equation}
 for some $\Phi(L_1^\pm, \cdots,
L_{n-2}^\pm)\in F[L_1^\pm, L_2^\pm, \cdots, L_{n-2}^\pm]\cap
\cba{n-2}$. On the other hand, by (\ref{center1}) , for any positive
integer $k$, we have
$$ E_{n-1}T_{n-2}L_{n-1}^k E_{n-1}=rL_{n-2}^k E_{n-1}+\omega\sum_{j=1}^k
(L_{n-2}^{j-1}\omega_{n-1}^{(k-j+1)}-L_{n-2}^{2j-k-2})
 E_{n-1}.
$$

Acting  $\sigma$ to  $E_{n-1}T_{n-2}L_{n-1}^k E_{n-1}$ yields the
formula for $E_{n-1}T^{-1}_{n-2}L_{n-1}^{-k} E_{n-1}$. Using
Lemma~\ref{wenzl-rel}(c) for $i=n-2$ to rewrite
$E_{n-1}T^{-1}_{n-2}L_{n-1}^{-k} E_{n-1}$ yields  the formula for
$E_{n-1} T_{n-2}L_{n-1}^{-k} E_{n-1}$. If $k=0$, then
$E_{n-1}T_{n-2}L_{n-1}^k E_{n-1}=r E_{n-1}$. So,  there is a $
\Psi(L_1^\pm,\cdots, L_{n-2}^\pm)\in F[L_1^\pm, L_2^\pm, \cdots,
L_{n-2}^\pm]\cap \cba{n-2}$ such that
\begin{equation}\label{Psi}
E_{n-1}T_{n-2}  F_{\t,n} F_{\t, n-1} E_{n-1}=\Psi (L_1^\pm,\cdots,
L_{n-2}^\pm) E_{n-1}.\end{equation} By Definition~\ref{def-m},
$E_{2f-1}E_{2f} \cdots E_{n-1} \widetilde\M_\lambda=\M_\u$ where
$\u\simnn \t$ with $\u_{n-1}=\lambda\cup \{(k+1,1)\}$ if $\mu_k>1$
and $\u_{n-1}=\t_{n-1}=\mu$ if $\mu_k = 1$. In the latter case,
$\u=\t$.

Let $\Phi_\lambda$ (resp. $\Psi_\lambda$)  be obtained  from $\Phi$
(resp. $\Psi$) by using $c_{\t^{\lambda}}^\pm(k)$ instead of
$L_k^\pm$ in $\Phi$ (resp.$\Psi$). Note that  $\widetilde \M_\lambda
T_j=q\widetilde \M_\lambda$ for $a\le j\le n-3$.  By
Theorem~\ref{xproduct} and the definition of $\u$,
$$
f_\t E_{n-1}= (\Phi_\lambda+q[\mu_k-1]\Psi_\lambda)\bar \M_\u
+\sum_{\v\overset {n-2}\succ \u } b_\v\bar \M_\v \prod_{k=1}^{n-2}
F_{\t, k}. $$ We use Lemma~\ref{fxproduct}(b) to express $\bar\M_\v$
as an $F$-linear combination of $f_\s$'s. So, $\s\succeq \v$. Note
that $L_i^{\pm}$ act on $f_\s$ as scalars for all $1\le i\le n$. So
is $\prod_{k=1}^{n-2} F_{\t, k} $. Since we are assuming that
$\v\overset{n-2} \succ \u $, $f_\u$ can not appear in the expression
of $\bar\M_\v \prod_{k=1}^{n-2} F_{\t, k} $. By
Lemma~\ref{fxproduct}(b), the coefficient of $f_\u$ in $\bar \M_\u$
is $1$. So, $E_{\t\u}(n-1)=\Phi_\lambda+q[\mu_k-1]\Psi_\lambda.$ We
assume $\t\neq \u$. Then
$$\begin{aligned} E_{\t\u} (n-1) f_\t E_{n-1} =  & (\Phi_\lambda+q[\mu_k-1]\Psi_\lambda) f_\t
E_{n-1}\\
=& f_\t E_{n-1}(1+q[\mu_k-1]T_{n-2}) F_{\t, n} F_{\t, n-1} E_{n-1}\quad\text{  by (\ref{Phi})-(\ref{Psi})} \\
= & \sum_{\v\simnn\t}
E_{\t\v}(n-1)f_\v(1+q[\mu_k-1]T_{n-2})F_{\t,n-1}F_{\t,n}E_{n-1}\\=
& E_{\t\t}(n-1) f_\t E_{n-1} +q^2[\mu_k-1] E_{\t\t}(n-1) f_\t
E_{n-1}\\
\end{aligned}
$$
In the last equality, we use $f_\v F_{\t, n}F_{\t, n-1}=0$ (resp.
$f_\v T_{n-2} F_{\t, n}F_{\t, n-1}=0$ )  for $\v\simnn \t$ and
$\v\neq\t$ which follows from Lemma~\ref{fxproduct}(d) (resp.
Lemma~\ref{fxproduct}(d) and Lemma~\ref{epro}(a)). We also use
Lemma~\ref{samerow}(a) to get $f_\t T_{n-2}=q f_\t$. So,
$E_{\t\u}(n-1)=[\mu_k]E_{\t\t}(n-1)$. We remark that the above
equality holds true when $\u=\t$. One can verify it similarly.  In
this case, $\mu_k=1$. Similar computation shows that
$\Phi_\lambda=E_{\t\t}(n-1)$ and $\Psi_\lambda= q E_{\t\t}(n-1)$.

By  similar arguments as above, we have $$\begin{aligned} &
f_{\t^{\lambda}\u}f_{\u \t^{\lambda}}\\ \equiv &
F_{\t^{\lambda}}E_{2f-1}T_{n, 2f}^{-1} T_{n-1, 2f-1}^{-1}
\widetilde\M_{\lambda} F_{\u,n-1}F_{\u,n}\widetilde
\M_{\lambda}T_{2f-1,n-1}^{-1} T_{2f, n}^{-1}
E_{2f-1}F_{\t^{\lambda}} \mod \cba{n}^{\vartriangleright (f,
\lambda)}\\  \equiv & F_{\t^{\lambda}} E_{2f-1} \cdots E_{n-2}
\widetilde\M_{\lambda} E_{n-1}F_{\u,n-1}F_{\u,n}E_{n-1}\widetilde
\M_{\lambda}
E_{n-2} \cdots E_{2f-1}F_{\t^{\lambda}} \mod \cba{n}^{\vartriangleright(f,\lambda)}\\
\equiv & E_{\u\u}(n-1) \delta^{f-1} [\lambda] ! F_{\t^\lambda}
E_{2f-1}\cdots E_{n-1}\widetilde \M_\lambda E_{n-2}\cdots E_{2f-1}
F_{\t^\lambda}\mod \cba{n}^{\vartriangleright(f, \lambda)}\\
= & E_{\u\u}(n-1) \delta^{f-1} [\lambda] ! F_{\t^\lambda} \bar
\M_\lambda T_{2f, n}T_{2f-1, n-1} E_{n-2}\cdots E_{2f-1}
F_{\t^\lambda}\\
= & E_{\u\u}(n-1) \delta^{f-1} [\lambda] ! F_{\t^\lambda} \bar
\M_\lambda
F_{\t^\lambda}\\
\equiv & E_{\u\u}(n-1)\delta^{f-1}[\lambda]!
f_{\t^{\lambda}{\t^\lambda}} \mod \cba{n}^{\vartriangleright(f,\lambda)}.\\
\end{aligned}
$$
 So, $\langle f_\u,
f_\u\rangle = E_{\u\u}(n-1) \delta^{f-1} [\lambda] !$. Since
$\langle \ , \ \rangle$ is associative, $\langle f_\u E_{n-1},
f_\t\rangle = \langle f_\u , f_\t E_{n-1}\rangle$. Thus $\langle
f_\u, f_\u\rangle E_{\t\u}(n-1) = \langle f_\t, f_\t\rangle
E_{\u\t}(n-1)$. By Lemma~\ref{ek}(b),
$$ \frac{\langle f_\t, f_\t\rangle }{
\delta^{f-1}[\mu]!}=  \frac{ 1 }{\delta^{f-1}[\mu]!}\frac{
[\mu_k]^2 E_{\t\t}(n-1)}{ E_{\u\u}(n-1)}E_{\u\u}(n-1) \delta^{f-1}
[\lambda] ! =[\mu_k] E_{\t\t}(n-1),$$ where $E_{\t\t}(n-1)$ can be
computed explicitly by Lemma~\ref{ek}(a).
\end{proof}

\begin{Prop}\label{key3}  Suppose $\t \in \UPD_n(\lambda)$ with
$(f, \lambda)\in \Lambda_n$, and $l(\lambda)=l$. If $\hat
\t=\t^{\mu}$, and $\t_{n-1} = \t_n \cup p$ with $p=(k,\mu_k)$
$k<l$,  define $\u= \t s_{n,a+1}$ with  $a = 2(f-1)+\sum_{j=1}^{k}
\mu_j$ and $\v = (\u_1,\cdots, \u_{a+1})$. Then
\begin{equation} \label{key4}
\frac{\langle f_\t,f_\t\rangle }{ \delta^{f-1}[\mu]!} =
\frac{[\mu_k] E_{\v\v}(a)} { r^2 q^{2(\mu_k-2k)}-1} \frac{
\prod_{r_1\in \mathscr A(\mu)^{<p}}( r^2 q^{-2(
c_{\mu}(p)-c_{\mu}(r_1))}-1)}{\prod_{r_2\in \mathscr R(\mu)^{<p}}
(r^2 q^{-2(c_{\mu}(p)+c_{\mu}(r_2))}-1)}.
\end{equation}
\end{Prop}

\begin{proof} By definition, $\u=\t s_{n, a+1}$.
Using the argument in the proof of Proposition~\ref{key0}, we have
$$ \begin{aligned} \langle f_\t, f_\t \rangle & = \langle f_\u,
f_\u \rangle  \prod_{j=a+2}^n \left(1-\omega^2
\frac{ c_\u(j)c_\u(a+1)}{(c_\u(j)-c_\u(a+1))^2}\right)\\
&=\frac {\langle f_\u, f_\u \rangle}{r^2 q^{2(\mu_k-2k)}-1}
 \frac {\prod_{r_1\in
\mathscr A(\mu)^{<p}} (r^2 q^{-2( c_{\mu}(p)-c_{\mu}(r_1))}-1)}{
\prod_{r_2\in
\mathscr R(\mu)^{<p}} ( r^2q^{-2( c_{\mu}(p)+c_{\mu}(r_2))}-1)}.\\
\end{aligned}
$$
 By Proposition~\ref{recursive}, $\langle f_\u,
f_\u \rangle=\langle f_{\v}, f_{\v} \rangle \prod_{i=k+1}^\ell
[\lambda_i]!$ where $\v\in \UPD_{a+1}(\nu)$ with $\v=(\u_1, \u_2,
\dots, \u_{a+1})$ and $\u_{a+1}=\nu$. Finally, we use
Proposition~\ref{key21} and $[\mu_k]!=[\lambda_k]![\mu_k]$ to get $
\langle f_{\v}, f_\v \rangle = E_{\v\v}(a) [\mu_k]^2 \delta^{f-1}
\prod_{i=1}^k [\lambda_i]!$. Simplifying $\langle f_\t, f_\t
\rangle$ via previous formulae yields (\ref{key4}), as required.
\end{proof}

\begin{Defn} Suppose $(f, \lambda)\in \Lambda_n$ and $(\ell, \mu)\in \Lambda_{n-1}$. Write
$(\ell, \mu)\rightarrow (f,\lambda)$ if either $\ell=f$ and
$\mu=\lambda\setminus \{p\}$ or $\ell=f-1$ and
$\mu=\lambda\cup\{p\}$. If $(\ell, \mu)\rightarrow (f, \lambda)$ we
define $\gamma_{\lambda/\mu}\in F$ to be the scalar by declaring
that
\begin{equation}\label{lammu} \gamma_{\lambda/\mu}= \frac{\langle f_\t,
f_\t\rangle}{\delta^\ell[\mu]!} \end{equation} where $\t\in
\UPD_n(\lambda)$ with $\hat \t=\t^\mu\in \UPD_{n-1}(\mu)$.
\end{Defn}

The following is the  first  main result of this paper.
\begin{Theorem}\label{main} Let $\cba{n}$ be the Birman-Wenzl algebra
over $\mathbb Z[r^\pm, q^\pm, \omega^{-1}]$, where $r, q$ are
indeterminates and $\omega=q-q^{-1}$. The Gram determinant $\det
G_{f, \lambda}$ associated to the cell module $\Delta(f, \lambda)$
of $\cba{n}$ can be computed by the following formula
\begin{equation}\label{maintheorem}\det G_{f, \lambda}=\prod_{(\ell, \mu)\rightarrow (f,
\lambda)} \det G_{\ell, \mu} \cdot \gamma_{\lambda/\mu}^{\dim
\Delta(\ell, \mu)}\in \mathbb Z[r^\pm, q^\pm,
\omega^{-1}].\end{equation} Furthermore, each scalar
$\gamma_{\lambda/\mu}$ can be computed explicitly  by (\ref{key1}),
(\ref{key2}), (\ref{key4}) and Lemma~\ref{ek}(a).
\end{Theorem}

\begin{proof}  We first compute the Gram
determinants over $\mathbb C(q^\pm, r^\pm, \omega^{-1})$. In order
to use the results in section~4, we have to use the fundamental
theorem of algebra (see \cite{RS:gram}).

Since $\tilde G_{f, \lambda}$, defined via orthogonal basis of
$\Delta(f, \lambda)$, is a diagonal matrix and  each diagonal is of
form $\langle f_\t, f_\t \rangle$, $\t\in \UPD_n(\lambda)$, we have
$\det \tilde G_{f, \lambda}=\prod_{\t\in \UPD_n(\lambda)} \langle
f_\t, f_\t \rangle$. By Proposition~\ref{recursive},
$$\det \tilde G_{f, \lambda}=\prod_{(\ell, \mu)\rightarrow (f,
\lambda)} \det\tilde  G_{\ell, \mu} \cdot \gamma_{\lambda/\mu}^{\dim
\Delta(l, \mu)}.$$ However, by Lemma~\ref{fxproduct}(f) $\det G_{f,
\lambda}=\det \tilde G_{f, \lambda}$ and $ \det  G_{\ell, \mu} =
\det\tilde  G_{\ell, \mu}$. Since the Jucys-Murphy basis of
$\Delta(f, \lambda)$ is defined over $\mathbb Z[r^\pm, q^\pm,
\omega^{-1} ]$, $\det G_{f, \lambda}\in \mathbb Z[r^\pm, q^\pm,
\omega^{-1}]$.
\end{proof}

Recall that $\lambda'$ is the dual partition of the partition
$\lambda$. The following result gives the relation between the
integral factors of $\det G_{f, \lambda}$ and   $\det G_{f,
\lambda'}$.

\begin{Cor}\label{symmetry} Let $\cba{n}$ be the Birman-Wenzl algebra
over $\mathbb Z[r^\pm, q^\pm, \omega^{-1}]$, where
$\omega=q-q^{-1}$. Suppose $(f, \lambda)\in \Lambda_n$ and
$\epsilon\in \{-1, 1\}$. Then $r-\epsilon q^a$ is a factor of $\det
G_{f, \lambda}$ if and only if $r+\epsilon q^{-a}$ is a factor of
$\det G_{f, \lambda'}$.
\end{Cor}

\begin{proof} Note that $p=(i, j)$  is an addable (resp. removable) node of $\lambda$ if and
only if $p'=(j, i)$ is an addable (resp. removable) node of
$\lambda'$.  Detailed analysis for the numerators and denominators
of $\gamma_{\lambda/\mu}$ by elementary computation yields the
result.\end{proof}

If we consider $q$ as a scalar, then  $\prod_{(\ell, \mu)\rightarrow
(f, \lambda)} \gamma_{\lambda/\mu}^{\dim \Delta(\ell, \mu)}$ can be
considered as a  rational function in $r$. Write
$\prod_{(\ell,\mu)\rightarrow (f, \lambda)}
\gamma_{\lambda/\mu}^{\dim \Delta(\ell, \mu)}=\frac{f(r)}{g(r)}$
such that the g.c.d of $f(r)$ and $g(r)$ is $1$. The following
result is useful when we determine the zero divisors of  Gram
determinants.

\begin{Cor}\label{integer} Let $f(r)$ and $g(r)$ be defined as above.  Suppose $\epsilon\in \{-1, 1\}$.
 If $r-\epsilon q^a\mid f(r)$ and
$r-\epsilon q^a\nmid g(r)$ with  $a\in \mathbb Z$, then
$r-\epsilon q^a$  is a factor of $\det G_{f, \lambda}$. In other
words, $\det G_{f, \lambda}=0$ if $r=\epsilon q^a$.\end{Cor}
\begin{proof} The result follows from the fact that $\det G_{f, \lambda}
\in  \mathbb Z[r^\pm, q^\pm, \omega^{-1}] $.
\end{proof}

\section{Semisimplicity  criteria for $\cba{n}$  over a field}
In this section, we consider  $\mathscr B_{n, F} $ over an arbitrary
field $F$. We will give a  necessary and sufficient condition for
$\mathscr B_{n, F}$  being (split) semisimple.  We will denote
$\mathscr B_{n, F}$ by $\cba{n}$ if there is no confusion.

We remark that  we may not have the orthogonal representations over
$F$. However, we still have the recursive formula in
(\ref{maintheorem}) since the Gram matrix associated to each cell
module is a matrix over $R$. In what follows, we will use this fact
frequently.

\begin{Prop}\label{semikey1} Let $\cba{n}, n\ge 2$ be the
Birman-Wenzl algebra over  $\mathbb Z[r^\pm, q^\pm, \omega^{-1}]$.
Then
\begin{equation}\label{line}\begin{aligned} \det G_{1, (n-2)}
 = & q^{\frac{1} {2}(n-1)(3n-4)}\left(\frac{[n-2]!}{r(q^2-1)}\right)^{\frac1 2 n(n-1)}
(r-q)^{\frac1 2 n(n-3)}\\
& \times (r+q^3)^{\frac12 (n-1)(n-2)}
 (r^2 -q^{6-2n})^{n-1}(r-q^{3-2n}).\\
 \end{aligned}
\end{equation}
\end{Prop}
\begin{proof} Let
\begin{itemize}\item $(\s_{k, 2})_i=(i)$ for $i\le k-1$ and $
(\s_{k, 2})_i=(i-2)$  for $k\le i\le n$.
\item Suppose $3\le j\le k$. Define  $(\s_{k, j})_i=(i)$ for $i\le j-2$ and $
(\s_{k, j})_{i}=(i-1, 1)$, for $j-1\le i\le k-1$, and $(\s_{k,
j})_i=(i-2)$, $k\le i\le n$.
\end{itemize}
 Then $\UPD_n(\lambda)=\{ \s_{k, j}\mid  2\le j\le k\le n\}$.
We use  Proposition~\ref{recursive} and (\ref{key1}), (\ref{key2}),
(\ref{key4}) to compute $\langle f_{\s_{k, j}}, f_{\s_{k,
j}}\rangle$. We have \begin{itemize}\item $\langle f_{\s_{2, 2}},
f_{\s_{2, 2}}\rangle = \delta [n-2]!$, \item $ \langle f_{\s_{k,
2}}, f_{\s_{k, 2}}\rangle =\frac{q^3[n-2]!}{r (q^2-1)}
\frac{r-q^{3-2k}}{r-q^{5-2k}} (r^2 q^{2k-6}-1)[k-1]$ for $3\le k\le
n$, \item $\langle f_{\s_{k,j}}, f_{\s_{k,j}}\rangle = \frac{q}{r}
\frac{[n-2]!}{q^2-1} (r-q)(r+q^3)\frac{[j-2]}{[j-1]}
\frac{r^2-q^{6-2k}}{r^2-q^{8-2k}}$ for  $3\le j\le k\le
n$.\end{itemize}
Thus
 \begin{equation}\label{2}
 \prod_{k=3}^n
\langle f_{\s_{k, 2}}, f_{\s_{k, 2}}\rangle
=\left(\frac{q^3[n-2]!}{r (q^2-1)}\right)^{n-2} [n-1]!
\frac{r-q^{3-2n}}{r-q^{-1}} \prod_{k=3}^n (r^2 q^{2k-6}-1),
\end{equation}
 and
\begin{equation}\label{3}\begin{aligned}
 \prod_{3\le j\le k\le n} \langle f_{\s_{k, j}}, f_{\s_{k,
j}}\rangle  = & \left(\frac{q}{r} \frac{[n-2]!}{q^2-1}
(r-q)(r+q^3)\right)^{\frac{(n-2)(n-1)}2}\\ & \times   \frac
{(r^2-q^{6-2n})^{n-2}}{[n-1]!} \prod_{k=3}^n \frac{1}
{r^2-q^{8-2k}}.\\ \end{aligned}
\end{equation}

Note that $\det G_{1, \lambda}=\prod_{2\le j\le k\le n}\langle
f_{s_{k, j}}, f_{s_{k,
 j}}\rangle$.  Now,  (\ref{line}) follows from elementary computation via equalities given above.
\end{proof}

Let $F$ be a field containing non-zero $\mathbf {q, r}$ and
$(\mathbf {q}-\mathbf {q}^{-1})^{-1}$.   The following results can
be verified directly.
\begin{Lemma}\label{sst} Suppose $o(\mathbf q^2)>n$ and
$(f, \lambda)\in \Lambda_n$.  For any $p_1\in \mathscr A(\lambda)$
and $p_2\in \mathscr R(\lambda)$,
$[c_\lambda(p_1)+c_\lambda(p_2)]\neq 0$.
\end{Lemma}

\begin{Prop} \label{semikey} Suppose that $n\ge 2$.   Let $\cba{n}$ be the Birman-Wenzl algebra
over a field $F$. Suppose $o(\mathbf {q}^2) > n$ and $ \mathbf{ r}
\not\in \{\mathbf{q}^{-1}, -\mathbf{q}\}$. Then $\cba{n}$ is
semisimple if and only if  $ \prod_{k=2}^n \det G_{1, (k-2)} \det
G_{1, (1^{k-2})} \ne 0$.
\end{Prop}

\begin{proof} We  use
Lemma~\ref{sst} to check that any factor $\q^a-\q^b$ in  the
numerators of (\ref{key1}), (\ref{key2}), (\ref{key4}) is not equal
to zero if $o(\q^2)>n$. We will  only consider the factors in $\det
G_{f, \lambda}$ with forms $\r\pm \q^a$, $a\in \mathbb Z$. When we
want to prove $\det G_{f, \lambda}=0$ for $\r=\epsilon \q^a$,
$\epsilon\in \{1, -1\}$, we will consider $\q$ to be the
indeterminate $q$
 first.  In other words, we have  $o(q)=\infty$. In order to
get the result for $n<o({\mathbf q^2})<\infty$, $\mathbf{q}\in F$,
we will specialize the indeterminate  $q$ to ${\q}\in F$. In the
remainder  of the proof, we consider $\q$ to be the  indeterminate
$q$. Also, we use $r$ instead of $\r$.

($\Rightarrow$) Suppose $ \det G_{1, (k-2)} = 0$ for some $2\le
k\le n$. We claim that there is a $(f, \lambda)\in \Lambda_n$ such
that $\det G_{f, \lambda} = 0$. It gives rise to a contradiction
since the Gram determinant associated to any cell module is   not
equal to zero if a cellular algebra is semisimple.

If $ \det G_{1, (1^{k-2})} = 0$, then we  use
Corollary~\ref{symmetry}  and the above claim twice to find a
partition $\lambda$ such that  $ \det G_{f, \lambda} =0 $. This will
give rise to  a contradiction, too. Consequently, $ \prod_{2\le k\le
n} \det G_{1, (k-2)} \det G_{1, (1^{k-2})} \ne 0$ if $\cba{n}$ is
semisimple.

Now, we prove our claim. Since $\cba{n}$ is semisimple,
 $\det G_{1, (n-2)} \neq 0$. It is easy to verify  $ \det G_{1, \emptyset} \ne 0$ if
  $ r \not\in\{q^{-1}, -q\}$. Therefore, we can assume $ 2<k <n$.
   By Proposition~\ref{semikey1}, $r \in \{ q^{3-2k},
\pm q^{3-k}, -q^3, q\}$, $3 < k < n$. The case $k=3$ has been dealt
with in \cite{Enyang} which says $r \in \{ q^{-3}, \pm 1 , -q^3\}$.
One can use   the  program  \cite{Gap} (written in GAP language) to
verify it directly. \footnote{F. Luebeck wrote the GAP
 program for Brauer algebras when he visited our department  in 2006.
 Imitating his program, we wrote the GAP program for Birman-Wenzl
 algebras. We take this opportunity to express our gratitude to him.}

\Case{1. $n - k$ is even} Let $f = \frac{n-k}{2} + 1$ and $\lambda =
(k - 2)$. Then $(f, \lambda)\in \Lambda_n$. If $(f-1,
\mu)\rightarrow (f, \lambda)$, then $\mu\in \{\mu_1, \mu_2\}$ and
$\mu_1 = (k - 1)$,   $\mu_2 = (k - 2, 1)$. By
Proposition~\ref{key21}, we have
 \begin{equation} \label{case1} r_{\lambda/\mu_1} r_{\lambda
/\mu_2}=\frac{q^{2k-2} [k-2] (r-q)(r+q^3)(r^2-q^{6-2k})^2
(r-q^{3-2k})} {r^2 (q^2-1)^2[k-1]
(r^2-q^{8-2k})(r-q^{5-2k})}\end{equation} When $k>3$,
$(r-q)(r+q^3)(r^2-q^{6-2k}) (r-q^{3-2k})$ and $(r^2-q^{8-2k})
(r-q^{5-2k})$  are co-prime\footnote{since  $o(q)=\infty$} in
$\mathbb Z[r^\pm, q^\pm, (q-q^{-1})^{-1}]$.  If $k=3$, then
$(r+q^3)(r^2-1) (r-q^{-3})$ and $(r+q) (r-q^{-1})$ are co-prime in
$\mathbb Z[r^\pm, q^\pm, (q-q^{-1})^{-1}]$.
 Consequently, by (\ref{maintheorem}),
 $\det
G_{(f-1, \mu_1)}\det G_{(f-1, \mu_2)}\det G_{f, (k-3)}$ has to be
divided by $(r^2-q^{8-2k})^{\dim \Delta(f-1, \mu_2)}
(r-q^{5-2k})^{\dim \Delta(f-1, \mu_1)}$. Therefore, $\det G_{f,
\lambda}=0$ if   $\det G_{1, (k-2)}=0$.

 \Case{2. $n - k$ is odd} There are three subcases we have to
discuss.

\subcase{2a.  $r = q$ for some $k > 3$ or $r = - q^3$ for some
integer $k$ with  $k
> 2$} Let $f=\frac{n-k+1}{2}$ and $ \lambda= (k-1)$. By (\ref{case1})
and Corollary~\ref{integer}, $\det G_{f, \lambda}=0$. We remark that
we can use (\ref{case1}) since   $\gamma_{\lambda/\mu}$ depends only
on $\lambda$ and $\mu$.

\subcase{2b. $r = q^{3-2k}$, for some integer $k$ with  $2<k<n$} Let
$f = \frac{n-k+1}{2}$ and $\lambda = ( k - 2, 1 )$. Suppose $k
> 3$. If  $(f-1, \mu) \rightarrow (f, \lambda)$, then $\mu\in
\{\mu_1,
 \mu_2, \mu_3\}$ where $\mu_1 = (k-1, 1)$, $\mu_2 = (k - 2, 2)$
 and $\mu_3 = (k - 2, 1, 1)$. By  Propositions~\ref{key21}, \ref{key3},
$$\begin{aligned}
 \gamma_{\lambda / \mu_1}& =
\frac{q^3(r^2q^{2k-4}-1)(r^2q^{2k-8}-1)(r-q^{3-2k})}{r(q^2-1)(r-q^{5-2k})(r^2q^{2k-6}-1)}
,\\ \gamma_{\lambda / \mu_2} &
=\frac{[k-3]q^2(rq-1)(r^2-q^{4-2k})(r^2-q^4)}{r[k-2](q^2-1)(r^2-q^{6-2k})(r-q)}
,\\ \gamma_{\lambda / \mu_3} &
=\frac{q[k-1](r+q^5)(r^2-q^{8-2k})(r^2-q^4)}{r[k][2](q^2-1)(r^2-q^{10-2k})(r+q^3)}.
\\ \end{aligned}$$

If $k = 3$, then  $\lambda = (1, 1)$. In this case, define  $\mu_1 =
(2, 1)$ and $\mu_2 = (1, 1, 1)$. By   Propositions~\ref{key21},
\ref{key3}, $$\gamma_{\lambda / \mu_1}
=\frac{q^2(r^2-q^2)(rq+1)}{r(r^2-1)(q^2-1)}(r-q^{-3}), \text{ and }\
\
 \gamma_{\lambda / \mu_2}
 =\frac{q (r+q^5)(r^2-q^2)}{r[3](q^2-1)(r+q^3)}. $$

\subcase{2c.  $r = \pm  q^{3-k}$, for some integer $k$ with $2 < k
< n$} If $k=3$, define $f = \frac{n}{2} $ and $\lambda =
\varnothing$. If $(f-1, \mu) \rightarrow (f, \lambda)$, then
$\mu=(1)$. If $(f-2, \nu) \rightarrow (f-1, \mu)$, then
$\nu\in\{\nu_1, \nu_2\}$ where
  $\nu_1= (2)$ and $\nu_2 = (1, 1)$.
By Propositions~\ref{key21}, \ref{key3},
$$ \gamma_{\lambda / \mu}  = \delta,\ \
  \gamma_{\mu / \nu_1}
=\frac{q^3(r-q^{-3})(r^2-1)}{r(q^2-1)(r-q^{-1})},\ \   \gamma_{\mu
/\nu_2} = \frac{q(r+q^3)(r^2-1)}{[2]r(q^2-1)(r+q)}.
$$
Note that $r\not\in \{q^{-1},  - q\}$. If $k = 4$, then  $r =
-q^{-1}$.  Let $\lambda=(1,1,1)$. If
 $(f-1, \mu) \to (f,\lambda)$, then $\mu\in \{\mu_1, \mu_2\}$
 where $\mu_1 = (2, 1, 1)$ and $\mu_2 = (1,1,1,1)$. By
 Propositions~\ref{key21}, \ref{key3},
 $$
r_{\lambda / \mu_1}= \frac{q^3 (r-q^{-3})(r^2-q^4) (r+q^{-1})}{r
(q^2-1)(r^2-q^2)}, \ \  r_{\lambda / \mu_2}  =
 \frac{q(r+q^{7}) (r^2-q^4)}{r[4](q^2-1)(r+q^5)}.
 $$

Suppose  $k > 4$. Let $f=\frac{n-k+1}{2}$ and $\lambda=(k-3, 1,
1)$. If $(f-1, \mu) \to (f, \lambda)$, then $\mu\in \{\mu_1,
\mu_2, \mu_3\}$ where
 $\mu_1 = (k - 2, 1, 1)$, $\mu_2 = (k -
3, 2, 1 )$ and $\mu_3 = (k - 3, 1, 1, 1)$. By
Propositions~\ref{key21}, \ref{key3},
$$\begin{aligned}
 \gamma_{\lambda / \mu_1}
&=\frac{q^{2k-5}(r-q^{5-2k})(r^2-q^{6-2k})(r^2-q^{12-2k})}
{r(q^2-1)(r^2-q^{10-2k})(r-q^{7-2k})}, \\
\gamma_{\lambda/ \mu_2}  &=
 \frac{[k-4]q^3(r-q^{-1}) (r^2-q^{6-2k})(r^2-q^6)(r+q)}{
r[k-3](q^2-1)(r^2-q^{8-2k})(r^2-q^4)},\\
  \gamma_{\lambda / \mu_3}
  & =\frac{q [k-1](r+q^7)(r^2-q^{12-2k}) (r^2-q^6)}
  {r[3][k](q^2-1)(r^2-q^{14-2k}) (r+q^5)}.\\
  \end{aligned}
  $$

In each case, by Corollary~\ref{integer}, $\det G_{f, \lambda}$ is
divided by either $(r-q^{3-2k})(r\pm q^{3-k})(r+q^3)(r-q)$ for
$k>3$ or $(r-q^{-3})(r\pm 1) (r+q^3)$ for $k=3$. Thus, $\det G_{f,
\lambda}=0$. This completes the proof of the claim.

($\Leftarrow$)  Suppose that $\cba{n}$ is not semisimple. Then $\det
G_{f, \lambda} = 0$ for some  $(f, \lambda) \in \Lambda_n$. By
(\ref{maintheorem}), either $\det G_{\ell, \mu}=0$ or the numerator
of $\gamma_{\lambda/\mu}$ is equal to zero for some $(\ell,
\mu)\rightarrow (f, \lambda)$. In the first case, by induction on
$n$ with $n\ge 3$, $\prod_{k=2}^{n-1} \det G_{1, (k-2)} \det G_{1,
(1^{k-2})}=0$, a contradiction. In the latter case, if  $l=f$, then
$\gamma_{\lambda/\mu}\neq 0$ since we are assuming that $o(q^2)>n$
and the numerator of  $\gamma_{\lambda/\mu }$ is a product of $[k]$
for $k\le n-1$. By  Propositions~\ref{key21}, \ref{key3}, we need
only consider the numerators of $\gamma_{\lambda/\mu}$ with
$\ell=f-1$. We claim that

a) $r^2q^{c_\lambda(p)+c_\lambda(q)}\neq 1$ for all  $p, q \in
\mathscr{A}(\lambda)$.

b) $r^2q^{c_{\mu}(q)-c_{\mu}(p)}\neq 1$ for all $q \in
\mathscr{A}(\mu)$ and $p\in \mathscr{R}(\mu)$.

At first, we prove a). We  assume that $p$ (resp. $q$) is in the
$k$-th (resp. $\ell$-th ) row. Since  two  boxes have the same
contents if they are in the same diagonal, we can move both $p$ and
$q$ to either the  first row or the first  column of a partition.
Therefore, there is a partition $\xi\in \{(k-2), (1^{k-2})\}$ for
some $2\leq k\leq n$ such that  $p_1, q_1\in \mathscr A(\xi)$ with
$c_{\lambda}(p)+c_\lambda(q)=c_{\xi}(p_1)+c_{\xi}(q_1)$. By our
assumption and Proposition~\ref{line} and Corollary~\ref{symmetry},
$c_{\lambda}(p)+c_{\lambda}(q)$ is a factor of $\det G_{1, \xi}$.
Since we are assuming that $\prod_{k=2}^n \det G_{1, (k-2)} \det
G_{1, (1^{k-2})} \ne 0$, we have  $\det G_{1, \xi}\neq 0$, forcing
$c_{\lambda}(p)+c_{\lambda}(q)\neq 0$.

b) can be proved by similar arguments as above.  We leave the
details to the reader. By Propositions~\ref{key21} and~\ref{key3},
and our claim b), $\gamma_{\lambda/\mu}=0$ implies either
$E_{\t\t}(n-1)=0$ or $E_{\v\v}(a)=0$, where $\t,  \v$ and $a$ are
defined in (\ref{lammu}) and  Proposition~\ref{key3}, respectively.

Using our claims a)-b), we have that $rq^{2c_\lambda(p_1)}-q^{-1}=0$
or $rq^{2c_\lambda(p_1)}+q=0$ if  $E_{\t\t}(n-1)E_{\v\v}(a)=0$. In
this case, $\mu$ is obtained from $\lambda$ by adding the addable
node $p_1$.

If $c_\lambda(p_1)=0$. then    $r\in \{q^{-1}, -q\}$, a
contradiction. So, we can assume that  $c_\lambda(p_1)\neq 0$.
First, we deal with the case when  $c_\lambda(p_1)>0$. Note that $ r
\in\{q^{-(1+2c_{\lambda}(p_1))}, -q^{1-2c_{\lambda}(p_1)}\}$. In the
first case, by (\ref{line}), $\det G_{1, \eta}=0$ where $\eta
=(c_\lambda(p_1))$. Note that $c_\lambda(p_1)\le n-2$, by
assumption, $\det G_{1, \eta}\neq 0$,  a contradiction. Assume that
$r= -q^{1-2c_{\lambda}(p_1)}$. By Proposition~\ref{line}, $\det
G_{1, \eta}=0$ where $\eta=(2c_\lambda(p_1))$. Since we are assuming
that  $\prod_{m=2}^n \det G_{1, (m-2)} \det G_{1, (1^{m-2})}\neq 0$,
we have $n-2<2c_\lambda(p_1)=2(\lambda_k+1-k)\le 2\lambda_k$,
forcing $k=1$. By (\ref{case1}), the numerators of
$\gamma_{\lambda/\mu_1} \gamma_{\lambda/\mu_2}$ must be divided by
$r+q^{1-2c_{\lambda}(p_1)}$, where $\mu_1, \mu_2$ is the same as
those in (\ref{case1}) and the $k$ in (\ref{case1}) should be
replaced by  $\lambda_1+2$ which is equal to $n-2f+2$.  Thus,
$r+q^{1-2c_{\lambda}(p_1)}\in S$ where
$$S= \{r-q, r+q^3, r\pm q^{1-\lambda_1},
r-q^{-(1+2\lambda_1)}\}.$$ On the other hand, we have
$\lambda_1+2\le n$ since we are assuming that   $f\ge 1$. By
(\ref{line}), $\prod_{m=2}^n \det G_{1, (m-2)} \det G_{1,
(1^{m-2})}$ is divided by  each element in $S$. This implies that
$r+q^{1-2c_{\lambda}(p_1)}\neq 0$, a contradiction.

If $c_\lambda(p_1)< 0$, we use Corollary~\ref{symmetry} to consider
$r\in \{-q^{1+2c_\lambda(p_1)}, q^{2c_\lambda(p_1)-1}\}$. In this
situation, we still get a contradiction by the  result  for
$c_\lambda(p_1)> 0$ stated above.
\end{proof}

\begin{Prop} \label{ssimple} Let $\cba{n}$
be the Birman-Wenzl algebra over a field $F$ containing the
parameters $\q^\pm, \r^\pm$ and $(\q-\q^{-1})^{-1}$. Assume  $ \r
\in \{\q^{-1}, -\q\}$.
\begin{enumerate}  \item
$\cba{n}$ is not semisimple if $n$ is either even or odd with
$n\ge 7$.
\item $\cba{1}$ is always semisimple.
\item $\cba{3}$ is semisimple if and only if $o(\q^2)>3$ and
$\q^4+1\neq 0$.
\item $\cba{5}$ is semisimple if and only if $o(\q^2)>5$ and
$\q^6+1\neq 0$, and $\q^8+1\neq 0$, and  $\Char F\neq 2$.
\end{enumerate}
\end{Prop}

\begin{proof} We  use $r, q$ instead of $\r, \q$
in the proof of this result. Since $r\in \{q^{-1}, -q\}$,
$\delta=0$. Suppose that  $n$ is even. Let $a=\dim_F \Delta (n/2 -
1, (1))$.
 By (\ref{maintheorem}),
$\det G_{n/2, \varnothing} = \det G_{n/2- 1, (1)} \delta^a=0$. We
have $\det G_{1, (3,2)}=0$ when $r\in \{q^{-1}, -q\}$.  One can use
\cite{Gap} to verify the  above  formulae easily. This shows that
$\cba{7}$ is not semisimple. We also use \cite{Gap} to get the
following formulae:
\begin{itemize} \item $
\det G_{1, (1)}=(q^4+1)$ if  $r\in \{q^{-1}, -q\}$.\item $\det G_{1,
(3)}= 2^5[2]^{10}[3]^{14}(1+q^8)$ (resp.
$-[2]^{10}[3]^{11}q^{-2}(1+q^4)^6$) if  $ r = -q$ (resp. if $ r =
q^{-1}$). \item  $\det G_{1, (1,1,1)}=q^{-2}[3](1+q^4)^6$ ( resp.
$2^5 [3]^4 (1+q^8)$) if  $ r = -q$ (resp. $ r = q^{-1}$).\item $\det
G_{1,(2, 1)}=-q^2[2]^4[3]^{15}
 (1+q^6)^4$ if $r\in \{q^{-1}, -q\}$.
\item
$\det G_{2, (1)}=-32 q^2(1+q^2) (1+q^4)^{10}(1+q^6)$ if $r\in
\{q^{-1}, -q\}$.
\end{itemize}
Now, (b)-(d) follow from the results on the semisimplicity of  Hecke
algebras $\H_n$ for $n\in \{1, 3, 5\}$ together with the above
formulae.

We close the proof by showing that  $\det G_{\frac{n-5}{2}, (3,
2)}=0$ for all odd $n$ with $n> 7$. This can be verified by
comparing the recursive formulae on $\det G_{\frac{n-5}{2}, (3, 2)}$
with $\det G_{1, (3, 2)}$. We leave the details to the reader.
\end{proof}

\begin{Theorem}\label{main2} Let $\cba{n}$ be the Birman-Wenzl algebra
over a field $F$ which contains non-zero  parameters $r, q, \omega$,
where  $\omega =q-q^{-1}$.
\begin{enumerate}
\item Suppose $ r \not\in \{q^{-1}, -q\}$.
\begin{itemize} \item [(a1)]If $n\ge 3$, then $\cba{n}$ is semisimple if and only if $o(q^2)
> n$ and $ r \not\in \cup_{k=3}^n \{q^{3-2k}, \pm q^{3-k},
-q^{2k-3}, \pm q^{k-3}\}$.
\item [(a2)] $\cba{2}$ is semisimple if and only if $o(q^2)>2$.
\item [(a3)] $\cba{1}$ is always semisimple.
\end{itemize}
\item Assume
$ r \in \{q^{-1}, -q\}.$ \begin{itemize} \item[(b1)] $\cba{n}$ is
not semisimple if  $n$ is either  even or odd with $n\ge 7$.
\item [(b2)] $\cba{1}$ is always semisimple.
\item [(b3)] $\cba{3}$ is semisimple if and only if $o(q^2)>3$ and
$q^4+1\neq 0$.
\item [(b4)] $\cba{5}$ is semisimple if and only if $o(q^2)>5$,
$q^6+1\neq 0$, and $q^8+1\neq 0$ and $\Char F\neq 2$.
\end{itemize}
\end{enumerate}
\end{Theorem}

\begin{proof}  Suppose $n\neq 2$. Theorem~\ref{main2} follows
from Propositions~\ref{semikey}, (\ref{line}) and
Corollary~\ref{symmetry}, (resp. Proposition~\ref{ssimple}) under
the assumption $r\not\in \{q^{-1}, -q\}$ (resp. $r\in \{q^{-1},
-q\}$).  When $n=2$, we compute $\det G_{1, \varnothing}$ directly
to verify the result.
\end{proof}

Let $\delta=\frac{(q+r)(qr-1)}{r(q+1)(q-1)}$. Then
$$\text{limit}_{q\to 1} \delta\in\{1, 2, \cdots, n-2\}\cup\{-2,
-4,\cdots, 4-2n\}\cup\{-1,-2, \cdots, 4-n\}$$ if $r\in \cup_{k=3}^n
\{q^{3-2k}, \pm q^{3-k}, -q^{2k-3}, \pm q^{k-3}\}$ and $n\ge 3$.
They are the parameters we got in \cite{RS:ssbrauer} such that the
corresponding Brauer algebra is not semisimple. Finally, we remark
that some partial results on Brauer algebras being semisimple over
$\mathbb C$ can be found in \cite{Brown:brauer, DoranHanlonWales,
Wenzl:ssbrauer}.

\providecommand{\bysame}{\leavevmode ---\ }
\providecommand{\og}{``} \providecommand{\fg}{''}
\providecommand{\smfandname}{and}
\providecommand{\smfedsname}{\'eds.}
\providecommand{\smfedname}{\'ed.}
\providecommand{\smfmastersthesisname}{M\'emoire}
\providecommand{\smfphdthesisname}{Th\`ese}

\end{document}